\documentstyle{article}

\newcommand{\eL}{\mbox{$\cal L$}}
\newcommand{\Pe}{\mbox{$\cal P$}}

\newcommand{\kon}{\wedge}
\newcommand{\dis}{\vee}
\newcommand{\str}{\rightarrow}
\newcommand{\rts}{\leftarrow}
\newcommand{\mj}{\mbox{\bf 1}}

\newcommand{\df}{\mbox{\scriptsize{\it df}}}
\newcommand{\HDS}{\vrule width0pt height2.3ex depth1.05ex\displaystyle}

\newcommand{\eM}{\mbox{$\cal M$}}

\newcommand{\mmds}{\mbox{\bf MDS}}

\newcommand{\mgds}{\mbox{\bf GDS}}
\newcommand{\mds}{\mbox{\bf DS}}

\def\cirk{\,{\raisebox{.3ex}{\tiny $\circ$}}\,}
\def\ks{\mbox{\footnotesize$\;\xi\;$}}
\def\b#1#2{\stackrel{\raisebox{-2pt}{\mbox{\tiny $#1$}}}
{\raisebox{0pt}{$b$}}^{\raisebox{-7pt}{\scriptsize $#2$}}}

\def\c#1{\stackrel{\raisebox{-2pt}{\mbox{\tiny $\,#1$}}}
{\raisebox{0pt}{$c$}}}

\def\f#1#2{{{\HDS #1}\over{\HDS #2}}}

\def\nav#1#2{\parshape=2 1em 31em 4em 28em \smallskip
\noindent{\makebox[3em][l]{#1}}{#2}\par\smallskip}

\def\prop#1#2{\vspace{2ex} \noindent{\sc #1.} {\it #2} \par \vspace{2ex}}
\def\dkz{\noindent{\sc Proof. }}
\def\qed{\hfill $\dashv$}

\def\pl{\!+\!}
\def\mn{\!-\!}

\def\n#1{n^{#1}}
\def\r#1#2{\stackrel{\raisebox{-2pt}{\mbox{\tiny $#1$}}}
{\raisebox{0pt}{$r$}}^{\raisebox{-3pt}{\scriptsize $#2$}}}

\def\kst{\raisebox{1pt}{\mbox{\tiny$\xi$}}}

\newcommand{\PN}{\mbox{\bf PN}}
\newcommand{\PNN}{$\PN^{\neg}$}

\newcommand{\MPN}{\mbox{\bf MPN}}
\newcommand{\MPNN}{$\MPN^{\neg}$}

\def\Del#1{\stackrel{\raisebox{-2pt}{\mbox{\tiny $#1$}}}
{\raisebox{0pt}{$\Delta$}}}
\def\Sig#1{\stackrel{\raisebox{-2pt}{\mbox{\tiny $#1$}}}
{\raisebox{0pt}{$\Sigma$}}}

\def\Dk{\Del{\kon}}
\def\Sd{\Sig{\dis}}
\def\Dkp{\Dk^{\raisebox{-7pt}{$'$}}}
\def\Sdp{\Sd^{\raisebox{-7pt}{$'$}}}

\def\Dd{\Del{\dis}}
\def\Dx{\Del{\xi}}
\def\Sk{\Sig{\kon}}
\def\Sx{\Sig{\xi}}

\def\Ddp{\Dd^{\raisebox{-7pt}{$'$}}}
\def\Dxp{\Dx^{\raisebox{-7pt}{$'$}}}
\def\Skp{\Sk^{\raisebox{-7pt}{$'$}}}
\def\Sxp{\Sx^{\raisebox{-7pt}{$'$}}}

\def\Ksi#1{\stackrel{\raisebox{-2pt}{\mbox{\tiny $#1$}}}
{\raisebox{0pt}{$\Xi$}}}
\def\Xk{\Ksi{\kon}}
\def\Xd{\Ksi{\dis}}
\def\Xx{\Ksi{\xi}}

\def\Xkp{\Xk^{\raisebox{-7pt}{$'$}}}
\def\Xdp{\Xd^{\raisebox{-7pt}{$'$}}}
\def\Xxp{\Xx^{\raisebox{-7pt}{$'$}}}

\def\The#1{\stackrel{\raisebox{-2pt}{\mbox{\tiny $#1$}}}
{\raisebox{0pt}{$\Theta$}}}
\def\Tk{\The{\kon}}
\def\Td{\The{\dis}}
\def \Fn{F^{\neg}}

\def\koc{\begin{picture}(10,7)
\put(1.5,0){\line(0,1){7}} \put(8.5,0){\line(0,1){7}}
\put(1.5,0){\line(1,0){7}} \put(1.5,7){\line(1,0){7}}
\end{picture}}

\def\koci{\begin{picture}(7,5)
\put(1,0){\line(0,1){5}} \put(6,0){\line(0,1){5}}
\put(1,0){\line(1,0){5}} \put(1,5){\line(1,0){5}}
\end{picture}}

\begin{document}

\title{Coherence of Proof-Net Categories}
\author{{\sc Kosta Do\v sen} and {\sc Zoran Petri\' c}\\[0.5cm]
Mathematical Institute, SANU \\
Knez Mihailova 35, p.f. 367 \\
11001 Belgrade, Serbia \\
email: \{kosta, zpetric\}@mi.sanu.ac.yu}
\date{}
\maketitle

\begin{abstract}
\noindent The notion of proof-net category defined in this paper
is closely related to graphs implicit in proof nets for the
multiplicative fragment without constant propositions of linear
logic. Analogous graphs occur in Kelly's and Mac Lane's coherence
theorem for symmetric monoidal closed categories. A coherence
theorem with respect to these graphs is proved for proof-net
categories. Such a coherence theorem is also proved in the
presence of arrows corresponding to the mix principle of linear
logic. The notion of proof-net category catches the unit free
fragment of the notion of star-autonomous category, a special kind
of symmetric monoidal closed category.
\end{abstract}

\vspace{.3cm}


\vspace{.5ex}

\noindent {\it Keywords$\,$}: generality of proofs, linear logic,
mix, proof nets, linear distribution, dissociativity, categorial
coherence, Kelly-Mac Lane graphs, Brauerian graphs, split
equivalences, symmetric monoidal closed category, star-autonomous
category

\vspace{1cm}

\baselineskip=1.2\baselineskip

\section{Introduction}

In this paper we introduce the notion of proof-net category, for
which we will show that it is closely related to graphs implicit
in proof nets for the multiplicative fragment without constant
propositions of linear logic (see \cite{G87} and \cite{DR89} for
the notion of proof net). Analogous graphs occur in Kelly's and
Mac lane's coherence theorem for symmetric monoidal closed
categories of \cite{KML71}.

The notion of proof-net category is based on the notion of
symmetric net category of \cite{DP04} (Section 7.6); these are
categories with two multiplications, $\kon$ and $\dis$,
associative and commutative up to isomorphism, which have moreover
arrows of the {\it dissociativity}\index{dissociativity} type
${A\kon(B\dis C)\str (A\kon B)\dis C}$ (called {\it linear} or
{\it weak} distribution\index{lzinear distribution@linear
distribution} in \cite{CS97}). The symmetric net category freely
generated by a set of objects is called \mds. To obtain proof-net
categories we add to symmetric net categories an operation on
objects corresponding to negation, which is involutive up to
isomorphism. With these operations come appropriate arrows. A
number of equations between arrows, of the kind called {\it
coherence conditions} in category theory, are satisfied in
proof-net categories.

A notion amounting to the notion of star-autonomous category of
\cite{Barr91} is obtained in a similar manner in \cite{CS97}.
Star-autonomous categories, which stem from \cite{Barr79}, are a
special kind of symmetric monoidal closed categories. In
contradistinction to symmetric net and proof-net categories they
involve unit objects.

We introduce next a category \emph{Br} whose arrows are called
{\it Brauerian split equivalences} of finite ordinals. These
equivalence relations, which stem from results in representation
theory of \cite{B37}, amount to the graphs used by
Kelly\index{Kelly@Kelly, G.M.} and Mac Lane\index{Mac Lane@Mac
Lane, S.} for their coherence theorem of symmetric monoidal
categories mentioned above. Brauerian split equivalences express
generality of proofs in linear logic (see \cite{DP03a},
\cite{DP03b}).

For proof-net categories we prove a coherence theorem that says
that there is a faithful functor from the proof-net category \PNN\
freely generated by a set of objects into \emph{Br}. The coherence
theorem for \PNN\ yields an elementary decision procedure for
verifying whether a diagram of arrows commutes in \PNN, and hence
also in every proof-net category. This is a very useful result,
which enables us in \cite{DP05} to obtain other coherence results
with respect to \emph{Br}, in particular a coherence result for
star-autonomous categories, involving the units. It is also shown
in \cite{DP05} with the help of coherence for \PNN\ that the
notion of proof-net category catches the unit-free fragment of
star-autonomous categories. (A different attempt to catch this
fragment is made in \cite{Lamar05} and \cite{HHS}.)

The coherence theorem for \PNN\ is proved by finding a category
\PN, equivalent to \PNN, in which negation can be applied only to
the generating objects, and coherence is first established for
\PN\ by relying on coherence for symmetric net categories,
previously established in \cite{DP04} (Chapter~7), and on an
additional normalization procedure involving negation.

In the last two sections of the paper we consider proof-net
categories that have \emph{mix}\index{mix} arrows of the type
${A\kon B\vdash A\dis B}$. We prove coherence with respect to
\emph{Br} for the appropriate notion of proof-net category with
these arrows, which we call mix-proof-net category.

\section{\large\bf The category \mds}

The objects of the category \mds\ are\index{DS@\mds\ category} the
formulae of the propositional language
${\eL_{\kon,\dis}}$,\index{Lcondis@${\eL_{\kon,\dis}}$ language}
generated from a set \Pe\ of propositional letters,\index{P@\Pe\
set of letters} which we call simply {\it
letters},\index{lzetter@letter} with the binary connectives
$\kon$\index{conjunction@conjunction $\kon$} and
$\dis$.\index{disjunction@disjunction $\dis$} We use
$p,q,r,\ldots\,$,\index{pb@$p$ letter}\index{q@$q$
letter}\index{r@$r$ letter} sometimes with indices, for letters,
and $A,B,C,\ldots\,$,\index{A@$A$ formula}\index{B@$B$
formula}\index{C@$C$ formula} sometimes with indices, for
formulae. As usual, we omit the outermost parentheses of formulae
and other expressions later on.

To define the arrows of \mds, we define first inductively a set of
expressions called the {\it arrow terms}\index{arrow term} of
\mds. Every arrow term of \mds\ will have a {\it
type},\index{type} which is an ordered pair of formulae of
${\eL_{\kon,\dis}}$. We write ${f\!:A\vdash
B}$\index{turnstile@turnstile $\vdash$} when the arrow term $f$ is
of type ${(A,B)}$. (We use the turnstile $\vdash$ instead of the
more usual $\str$, which we reserve for a connective and a
biendofunctor.) We use $f,g,h,\ldots\,$, sometimes with indices,
for arrow terms.\index{f@$f$ arrow}\index{G functor 1@$g$
arrow}\index{h@$h$ arrow}

For all formulae $A$, $B$ and $C$ of ${\eL_{\kon,\dis}}$ the
following {\it primitive arrow terms}:\index{primitive arrow
terms}

\begin{tabbing}
\centerline{$\mj_A\!: A\vdash A$,}\index{identity arrow@identity
arrow $\mj$}
\\[1ex]
\mbox{\hspace{.2em}}\= $\b{\kon}{\str}_{A,B,C}\,$\= : \=
$A\kon(B\kon C)\vdash (A\kon B)\kon C$,\quad \=
$\b{\dis}{\str}_{A,B,C}\,$\= : \= $A\dis(B\dis C)\vdash (A\dis
B)\dis C$,\kill \> $\b{\kon}{\str}_{A,B,C}\,$\> : \> $A\kon(B\kon
C)\vdash (A\kon B)\kon C$,\> $\b{\dis}{\str}_{A,B,C}\,$\> : \>
$A\dis(B\dis C)\vdash (A\dis B)\dis
C$,\index{bcona@$\b{\kon}{\str}$
arrow}\index{bdisa@$\b{\dis}{\str}$ arrow}
\index{betxia@$\b{\xi}{\str}$ arrow}
\\*[1ex]
\> $\b{\kon}{\rts}_{A,B,C}\,$\> : \> $(A\kon B)\kon C\vdash A\kon
(B\kon C)$,\> $\b{\dis}{\rts}_{A,B,C}\,$\> : \> $(A\dis B)\dis
C\vdash A\dis (B\dis C)$,\index{bconb@$\b{\kon}{\rts}$
arrow}\index{bdisb@$\b{\dis}{\rts}$ arrow}
\index{betxib@$\b{\xi}{\rts}$ arrow}
\\[1ex]
\> \> $\c{\kon}_{A,B}\!\!$\' : \> $A\kon B\vdash B\kon A$,\> \>
$\c{\dis}_{A,B}\!\!$\' : \> $B\dis A\vdash A\dis
B$,\index{ccon@$\c{\kon}$ arrow}\index{cdis@$\c{\dis}$ arrow}
\\[1ex]
\centerline{$d_{A,B,C}\!: A\kon(B\dis C)\vdash (A\kon B)\dis
C$}\index{d@$d$ arrow}
\end{tabbing}

\noindent are arrow terms of \mds. If ${g\!:A\vdash B}$ and
${f\!:B\vdash C}$ are arrow terms of \mds, then ${f\cirk
g\!:A\vdash C}$\index{composition@composition $\cirk$} is an arrow
term of \mds; and if ${f\!:A\vdash D}$ and ${g\!:B\vdash E}$ are
arrow terms of \mds, then ${f\ks g\!:A\ks B\vdash D\ks E}$, for
${\!\ks\!\in\{\kon,\dis\}}$,\index{xi@$\ks$ connective} is an
arrow term of \mds. This concludes the definition of the arrow
terms of \mds.

Next we define inductively the set of {\it equations} of \mds,
which are expressions of the form ${f=g}$, where $f$ and $g$ are
arrow terms of \mds\ of the same type. We stipulate first that all
instances of ${f=f}$ and of the following equations are equations
of \mds:

\begin{tabbing}
\mbox{\hspace{1em}}\= $({\mbox{{\it cat}~1}})$\quad\quad\= $f\cirk
\mj_A=\mj_B\cirk f=f\!:A\vdash B$,\index{cat1@({\mbox{{\it
cat}~1}}) equation}
\\*[1ex]
\> $({\mbox{{\it cat}~2}})$\> $h\cirk (g\cirk f)=(h\cirk g)\cirk
f$,\index{cat2@({\mbox{{\it cat}~2}}) equation}
\\[1.5ex]
for $\!\ks\!\in\{\kon,\dis\}$,
\\*[1ex]
 \> $(\!\ks\, 1)$\> $\mj_A\ks\mj_B=\mj_{A\kst B}$,\index{xi1@${(\ks\, 1)}$ equation}
 \index{conjunction1@${(\kon\, 1)}$ equation} \index{disjunction1@${(\dis\, 1)}$ equation}
\\*[1ex]
\> $(\!\ks\, 2)$\> $(g_1\cirk f_1)\ks(g_2\cirk f_2)=(g_1\ks
g_2)\cirk(f_1\ks f_2)$,\index{xi2@${(\ks\, 2)}$ equation}
 \index{conjunction2@${(\kon\, 2)}$ equation} \index{disjunction2@${(\dis\, 2)}$ equation}
\\[1.5ex]
for $f\!:A\vdash D$, $g\!:B\vdash E$ and $h\!:C\vdash F$,
\\*[1ex]
\> $\mbox{($\b{\xi}{\str}$ {\it nat})}$\>  $((f\ks g)\ks
h)\cirk\!\b{\xi}{\str}_{A,B,C}\:=\:
\b{\xi}{\str}_{D,E,F}\!\cirk(f\ks (g\ks
h))$,\index{betxib1@\mbox{($\b{\xi}{\str}$ {\it nat})} equation}
\index{bconb1@\mbox{($\b{\kon}{\str}$ {\it nat})} equation}
\index{bdisba1@\mbox{($\b{\dis}{\str}$ {\it nat})} equation}
\\[1ex]
\> $\mbox{($\c{\kon}$ {\it nat})}$\> $(g\kon
f)\cirk\!\c{\kon}_{A,B}\:=\:\c{\kon}_{D,E}\!\cirk(f\kon
g)$,\index{cconnat@\mbox{($\c{\kon}$ {\it nat})} equation}
\\[1ex]
\> $\mbox{($\c{\dis}$ {\it nat})}$\> $(g\dis
f)\cirk\!\c{\dis}_{B,A}\:=\:\c{\dis}_{E,D}\!\cirk(f\dis
g)$,\index{cdisnat@\mbox{($\c{\dis}$ {\it nat})} equation}
\\[1ex]
\> $\mbox{($d$ {\it nat})}$\> $((f\kon g)\dis h)\cirk d_{A,B,C} =
d_{D,E,F}\cirk(f\kon (g\dis h))$,\index{dassocnat@\mbox{($d$ {\it
nat})} equation}
\\[1.5ex]
\pushtabs $(\b{\xi}{}\b{\xi}{})$\quad\=
$\b{\xi}{\rts}_{A,B,C}\!\cirk\!\b{\xi}{\str}_{A,B,C}\;=\mj_{A\kst(B\kst
C)}$,\quad\quad\=
$\b{\xi}{\str}_{A,B,C}\!\cirk\!\b{\xi}{\rts}_{A,B,C}\;
=\mj_{(A\kst B)\kst C}$,\index{betxibxi@${(\b{\xi}{}\b{\xi}{})}$
equation} \index{bconb4@${(\b{\kon}{}\b{\kon}{})}$
equation}\index{bdisba4@${(\b{\dis}{}\b{\dis}{})}$ equation}
\\[1ex]
$(\b{\xi}{}\!5)$\> $\b{\xi}{\rts}_{A,B,C\kst D}\!\cirk\!
\b{\xi}{\rts}_{A\kst
B,C,D}\;=(\mj_A\:\ks\b{\xi}{\rts}_{B,C,D})\cirk\!\b{\xi}{\rts}_{A,B\kst
C,D}\!\cirk(\b{\xi}{\rts}_{A,B,C}\ks\:\mj_D)$,\index{betxib3@\mbox{$(\b{\xi}{}5)$}
equation} \index{bconb3@\mbox{$(\b{\kon}{}5)$}
equation}\index{bdisba3@\mbox{$(\b{\dis}{}5)$} equation}
\\[2ex]
$(\c{\kon}\c{\kon})$\>
$\c{\kon}_{B,A}\!\cirk\!\c{\kon}_{A,B}\;=\mj_{A\kon
B}$,\index{cconccon@${(\c{\kon}\c{\kon})}$ equation}
\\*[1ex]
$(\c{\dis}\c{\dis})$\quad
$\c{\dis}_{A,B}\!\cirk\!\c{\dis}_{B,A}\;=\mj_{A\dis
B}$,\index{cdiscdis@${(\c{\dis}\c{\dis})}$ equation}
\\[1.5ex]
$(\b{\kon}{}\c{\kon})$\> $(\mj_B\:\kon \c{\kon}_{C,A})\cirk\!
\b{\kon}{\rts}_{B,C,A}\!\cirk\!\c{\kon}_{A,B\kon
C}\!\cirk\!\b{\kon}{\rts}_{A,B,C}\!\cirk(\c{\kon}_{B,A}\kon\:\mj_C)=
\;\b{\kon}{\rts}_{B,A,C}$,\index{bconb5@${(\b{\kon}{}\c{\kon})}$
equation}
\\*[1ex]
$(\b{\dis}{}\c{\dis})$\> $(\mj_B\:\dis \c{\dis}_{A,C})\cirk\!
\b{\dis}{\rts}_{B,C,A}\!\cirk\!\c{\dis}_{B\dis C,A}
\!\cirk\!\b{\dis}{\rts}_{A,B,C}\!\cirk(\c{\dis}_{A,B}\dis\:\mj_C)=
\;\b{\dis}{\rts}_{B,A,C}$,\index{bdiscdis@${(\b{\dis}{}\c{\dis})}$
equation}
\\[1.5ex]
$(d \kon)$\> $(\b{\kon}{\rts}_{A,B,C}\dis\:\mj_D)\cirk d_{A\kon
B,C,D}=d_{A,B\kon C,D}\cirk(\mj_A\kon
d_{B,C,D})\cirk\!\b{\kon}{\rts}_{A,B,C\dis
D}$,\index{dassoccon@${(d \kon)}$ equation}
\\[1ex]
$(d \dis)$\> $d_{D,C,B\dis
A}\cirk(\mj_D\:\kon\b{\dis}{\rts}_{C,B,A})=\;\b{\dis}{\rts}_{D\kon
C,B,A}\!\cirk(d_{D,C,B}\dis\mj_A)\cirk d_{D,C\dis
B,A}$,\index{dassocdis@${(d \dis)}$ equation}
\\[1.5ex]
for $d^R_{C,B,A}=_{\df}\;\c{\dis}_{C,B\kon A} \!\cirk
(\c{\kon}_{A,B} \dis\:\mj_C)\cirk
d_{A,B,C}\cirk(\mj_A\:\kon\c{\dis}_{B,C})\cirk\!\c{\kon}_{C\dis
B,A}:$
\\*[.5ex]
\` $(C\dis B)\kon A\vdash C\dis(B\kon A)$,\index{dassocr@$d^R$
arrow}
\\[1ex]
$(d\!\b{\kon}{})$\> $d^R_{A\kon
B,C,D}\cirk(d_{A,B,C}\kon\mj_D)=d_{A,B,C\kon D}\cirk(\mj_A\kon
d^R_{B,C,D})\cirk\!\b{\kon}{\rts}_{A,B\dis
C,D}$,\index{dassocbcon@${(d\b{\kon}{})}$ equation}
\\*[1ex]
$(d\!\b{\dis}{})$\>$(\mj_D\dis d_{C,B,A})\cirk d^R_{D,C,B\dis
A}=\;\b{\dis}{\rts}_{D,C\kon
B,A}\!\cirk(d^R_{D,C,B}\dis\mj_A)\cirk d_{D\dis
C,B,A}$.\index{dassocbdis@${(d\b{\dis}{})}$ equation} \poptabs
\end{tabbing}

The set of equations of \mds\ is closed under symmetry and
transitivity of equality and under the rules

\[
(\mbox{\it cong~}\ks\!)\quad \f{f=f_1 \quad \quad \quad g=g_1}
{f\ks g=f_1\ks g_1}\index{congxi@(\mbox{\it cong~}\ks) rule}
\]

\noindent where ${\!\ks\!\in\{\cirk,\kon,\dis\}}$, and if
$\!\ks\!$ is $\cirk$, then ${f\cirk g}$ is defined (namely, $f$
and $g$ have appropriate, composable, types).

On the arrow terms of \mds\ we impose\index{impose equations} the
equations of \mds. This means that an arrow of \mds\ is an
equivalence class of arrow terms of \mds\ defined with respect to
the smallest equivalence relation such that the equations of \mds\
are satisfied (see \cite{DP04}, Section 2.3, for details).

The equations $(\!\ks\, 1)$ and $(\!\ks\, 2)$ are called {\it
bifunctorial}\index{bifunctorial equations} equations. They say
that $\kon$ and $\dis$ are biendofunctors (i.e.\ 2-endofunctors in
the terminology of \cite{DP04}, Section 2.4).\index{conjunction
biendofunctor@$\kon$ biendofunctor}\index{disjunction
biendofunctor@$\dis$ biendofunctor}

It is easy to show that for \mds\ we have the equations

\begin{tabbing}
\mbox{\hspace{2em}}\= $\mbox{($\b{\xi}{\rts}$ {\it nat})}$\quad\=
$(f\ks(g\ks
h))\cirk\!\b{\xi}{\rts}_{A,B,C}\;=\;\b{\xi}{\rts}_{D,E,F}\!\cirk((f\ks
g)\ks h)$,\index{betxib2@\mbox{($\b{\xi}{\rts}$ {\it nat})}
equation} \index{bconb2@\mbox{($\b{\kon}{\rts}$ {\it nat})}
equation} \index{bdisba2@\mbox{($\b{\dis}{\rts}$ {\it nat})}
equation}
\\[1ex]
\> $\mbox{($d^R$ {\it nat})}$\> $(h\dis(g\kon f))\cirk
d^R_{C,B,A}=d^R_{F,E,D}\cirk((h\dis g)\kon
f)$.\index{dassocrnat@\mbox{($d^R$ {\it nat})} equation}
\end{tabbing}

\noindent We call these equations and other equations with
``\emph{nat}'' in their names, like those in the list above, {\it
naturality}\index{nbtt@naturality equations} equations. Such
equations say that $\b{\kon}{\str}$, $\b{\kon}{\rts}$, $\c{\kon}$,
etc.\ are natural transformations.

The equations $(d \kon)$, $(d \dis)$, ${(d\!\b{\kon}{})}$ and
${(d\!\b{\dis}{})}$ stem from \cite{CS97} (Section 2.1; see
\cite{CS92}, Section 2.1, for an announcement). The equation
${(d\!\b{\dis}{})}$ of \cite{DP04} (Section 7.2) amounts with
${(\b{\dis}{}\b{\dis}{})}$ to the present one.

\section{\large\bf The category \PNN}

The category \PNN\ is\index{PNN@\PNN\ category} defined as \mds\
save that we make the following changes and additions. Instead of
${\eL_{\kon,\dis}}$, we have the propositional language
${\eL_{\neg,\kon,\dis}}$,\index{Lnegcondis@${\eL_{\neg,\kon,\dis}}$
language} which has in addition to what we have for
${\eL_{\kon,\dis}}$ the unary connective
$\neg$.\index{negation@negation $\neg$}

To define the arrow terms of \PNN, in the inductive definition we
had for the arrow terms of \mds\ we assume in addition that for
all formulae $A$ and $B$ of ${\eL_{\neg,\kon,\dis}}$ the following
{\it primitive arrow terms}\index{primitive arrow terms}:

\begin{tabbing}
\centerline{$\Dk_{B,A}:A\vdash A\kon(\neg B\dis
B)$,}\index{Deltacon@$\Dk$ arrow}
\\*[1ex]
\centerline{$\Sd_{B,A}:(B\kon\neg B)\dis A\vdash
A$,}\index{Sigmadis@$\Sd$ arrow}
\end{tabbing}

\noindent are arrow terms of \PNN. We call the index $B$, of
$\Dk_{B,A}$ and $\Sd_{B,A}$ the {\it crown}\index{crown index}
index, and $A$ the {\it stem}\index{stem index} index. The right
conjunct ${\neg B\dis B}$ in the target of ${\Dk_{B,A}:A\vdash
A\kon(\neg B\dis B)}$ is the {\it crown}\index{crown} of
${\Dk_{B,A}}$, and the left disjunct ${B\kon\neg B}$ in the source
of ${\Sd_{B,A}:(B\kon\neg B)\dis A\vdash A}$ is the {\it
crown}\index{crown} of $\Sd_{B,A}$. We have analogous definitions
of crown and stem indices, and crowns for $\Sk$, $\Dkp$, $\Skp$,
$\Dd$, $\Sdp$ and $\Ddp$, which will be defined below. (The symbol
$\Delta$ should be associated with the Latin {\emph dexter},
because in $\Dk_{B,A}$, $\Dkp_{B,A}$, $\Dd_{B,A}$ and $\Ddp_{B,A}$
the crown is on the right-hand side of the stem; analogously,
$\Sigma$ should be associated with {\emph sinister}.)

To define the arrows of \PNN, we assume in the inductive
definition we had for the equations of \mds\ the following
additional equations, which we call the \PN\ {\it
equations}\index{PN equations@\PN\ equations} (and {\emph not}
\PNN\ equations):

\begin{tabbing}
\mbox{\hspace{2em}}\= $\mbox{($\Dk$ {\it nat})}$\quad\=
$(f\kon\mj_{\neg B\dis B})\cirk\!\Dk_{B,A}\;=\;\Dk_{B,D}\!\cirk
f$,\index{Deltaconnat@\mbox{($\Dk$ {\it nat})} equation}
\\*[1ex]
\> $\mbox{($\Sd$ {\it nat})}$\> $f\cirk\!
\Sd_{B,A}\;=\;\Sd_{B,D}\!\cirk(\mj_{B\kon \neg B}\dis
f)$,\index{Sigmadisnat@\mbox{($\Sd$ {\it nat})} equation}
\\[2ex]
\> $\mbox{($\b{\kon}{}\Dk$)}$\> $\b{\kon}{\rts}_{A,B,\neg C\dis
C}\!\cirk\!\Dk_{C,A\kon
B}\;=\mj_A\:\kon\Dk_{C,B}$,\index{bcondeltacon@\mbox{($\b{\kon}{}\Dk$)}
equation}
\\*[1ex]
\> $\mbox{($\b{\dis}{}\Sd$)}$\> $\Sd_{C,B\dis
A}\!\cirk\!\b{\dis}{\rts}_{C\kon\neg
C,B,A}\;=\;\Sd_{C,B}\dis\:\mj_A$,\index{bdissigmadis@\mbox{($\b{\dis}{}\Sd$)}
equation}
\\[2ex]
\quad for $\Sk_{B,A}\;=_{\df}\;\:\c{\kon}_{A,\neg B\dis
B}\!\cirk\!\Dk_{B,A}\: : A\vdash (\neg B\dis B)\kon
A$,\index{Sigmacon@$\Sk$ arrow}
\\*[1.5ex]
\> $\mbox{($d\!\Sk$)}$\> $d_{\neg A\dis A,B,C}\cirk\!\Sk_{A,B\dis
C}\;=\;\Sk_{A,B}\dis\:\mj_C$,\index{dassocnat1@\mbox{($d\Sk$)}
equation}
\\[2ex]
\quad for $\Dd_{B,A}\;
=_{\df}\;\:\Sd_{B,A}\!\cirk\!\c{\dis}_{B\kon\neg B,A}\:
:A\dis(B\kon\neg B)\vdash A$,\index{Deltadis@$\Dd$ arrow}
\\*[1.5ex]
\> $\mbox{($d\!\Dd$)}$\> $\Dd_{A,C\kon B}\!\cirk d_{C,B,A\kon\neg
A}\;=\;\mj_C\:\kon\Dd_{A,B}$,\index{dassocdeltadis@\mbox{($d\Dd$)}
equation}
\\[2ex]
\> $\mbox{($\Sd\Dk$)}$\> $\Sd_{A,A}\!\cirk d_{A,\neg A,A}
\cirk\!\Dk_{A,A}\;=\;\mj_A$,\index{Sigmadisdeltacon@\mbox{($\Sd\Dk$)}
equation}
\\[2ex]
\quad for \=$\Dkp_{B,A}\;=_{\df}\;(\mj_A\:\kon\c{\dis}_{B,\neg
B})\cirk\!\Dk_{B,A}\: :A\vdash A\kon(B\dis\neg
B)$\index{Deltaconprime@$\Dkp$ arrow} and
\\*[1ex]
\>$\Sdp_{B,A}\;=_{\df}\;\:\Sd_{B,A}\!\cirk(\c{\kon}_{\neg
B,B}\dis\:\mj_A):(\neg B\kon B)\dis A\vdash
A$,\index{Sigmadisprime@$\Sdp$ arrow}
\\*[2ex]
\mbox{\hspace{2em}}$\mbox{($\Sdp\Dkp$)}$\>\> $\Sdp_{A,\neg
A}\!\cirk d_{\neg A,A,\neg A}\cirk\!\Dkp_{A,\neg A}\;=\mj_{\neg
A}$.\index{Sigmadisprimedeltaconprime@\mbox{($\Sdp\Dkp$)}
equation}
\end{tabbing}

It is easy to show that for \PNN\ we have the equations

\begin{tabbing}
\mbox{\hspace{2em}}\= $\mbox{($\b{\xi}{\rts}$ {\it nat})}$\quad\=
$(f\ks(g\ks
h))\cirk\!\b{\xi}{\rts}_{A,B,C}\;=\;\b{\xi}{\rts}_{D,E,F}\!\cirk((f\ks
g)\ks h)$,\kill

\> $\mbox{($\Sk$ {\it nat})}$\> $(\mj_{\neg B\dis B}\kon
f)\cirk\!\Sk_{B,A}\;=\;\Sk_{B,D}\!\cirk
f$,\index{Sigmaconnat@\mbox{($\Sk$ {\it nat})} equation}
\\[1ex]
\> $\mbox{($\Dd$ {\it nat})}$\>
$f\cirk\!\Dd_{B,A}\;=\;\Dd_{B,D}\!\cirk(f\dis\mj_{B\kon\neg
B})$.\index{Deltadisnat@\mbox{($\Dd$ {\it nat})} equation}
\end{tabbing}

\noindent The naturality equations \mbox{($\Dk$ {\it nat})} and
\mbox{($\Sd$ {\it nat})} together with these say that $\Dk$,
$\Sd$, $\Sk$ and $\Dd$ are natural transformations in the stem
index only, i.e.\ in the second index.

We also have the following abbreviations:

\begin{tabbing}
\quad\quad\quad\quad\quad\quad\quad\=$\Skp_{B,A}\;$\=$=_{\df}\;\:
\c{\kon}_{A,B\dis\neg B}\!\cirk\!\Dkp_{B,A}\:$\=$ :A\vdash
(B\dis\neg B)\kon A$,\index{Sigmaconprime@$\Skp$ arrow}
\\[1ex]
\>$\Ddp_{B,A}$\>$=_{\df}\;\; \Sdp_{B,A}\!\cirk\!\c{\dis}_{\neg
B\kon B,A}$\>$ :A\dis(\neg B\kon B)\vdash
A$.\index{Deltadisprime@$\Ddp$ arrow}
\end{tabbing}

\noindent If $\Xi$\index{Xi@$\Xi$ arrow} stands for either
$\Delta$ or $\Sigma$ and $\!\ks\!\in\{\kon,\dis\}$, then for every
$\mbox{($\Xx$ {\it nat})}$\index{Xixinat@$\mbox{($\Xx$ {\it
nat})}$ equation} equation we have in \PNN\ the equation
$\mbox{($\Xxp$ {\it nat})}$,\index{Xixiprimenat@$\mbox{($\Xxp$
{\it nat})}$ equation} which differs from $\mbox{($\Xx$ {\it
nat})}$ by replacing $\Xx$\index{Xixi@$\Xx$
arrow}\index{Xicon@$\Xk$ arrow}\index{Xidis@$\Xd$ arrow} by
$\Xxp$,\index{Xixiprime@$\Xxp$ arrow}\index{Xiconprime@$\Xkp$
arrow}\index{Xidisprime@$\Xdp$ arrow} and the index of $\mj$ by
the appropriate index. For example, we have

\begin{tabbing}
\mbox{\hspace{2em}}\= $\mbox{($\Dkp$ {\it nat})}$\quad\quad\=
$(f\kon\mj_{B\dis\neg B})\cirk\!\Dkp_{B,A}\;=\;\Dkp_{B,D}\!\cirk
f$.\index{Deltaconprimenat@$\mbox{($\Dkp$ {\it nat})}$ equation}
\end{tabbing}

\noindent As alternative primitive arrow terms for defining \PNN\
we could take one of $\Xk$ or $\Xkp$ and one of $\Xd$ or $\Xdp\!$.

We can also derive for \PNN\ the following equations:

\begin{tabbing}
\mbox{\hspace{2em}}\= $\mbox{($\Dkp$ {\it nat})}$\quad\quad\=
$(f\kon\mj_{B\dis\neg B})\cirk\!\Dkp_{B,A}\;=\;\Dkp_{B,D}\!\cirk
f$.\kill

\mbox{\hspace{2em}}\> $\mbox{($\b{\kon}{}\Dk\Sk$)}$\>
$\b{\kon}{\rts}_{A,\neg B\dis B,C}\!\cirk(\Dk_{B,A}\kon\:\mj_C)=
\mj_A\:\kon\Sk_{B,C}$,\index{bcondeltaconsigmacon@(${\b{\kon}{}\Dk\Sk}$)
equation}
\\[1ex]
\> $\mbox{($\b{\kon}{}\Sk$)}$\> $\b{\kon}{\str}_{\neg C\dis
C,B,A}\!\cirk\!\Sk_{C,B\kon
A}\;=\;\Sk_{C,B}\kon\:\mj_A$.\index{bconsigmacon@(${\b{\kon}{}\Sk}$)
equation}
\end{tabbing}

\noindent For the first equation, with indices omitted, we have

\begin{tabbing}
\mbox{\hspace{2em}}\= $\b{\kon}{\rts}\!\cirk(\Dk\kon\;\mj)\;$\=
$=\;\b{\kon}{\rts}\!\cirk\!\c{\kon}\!\cirk(\mj\:\kon\Dk)\cirk\!\c{\kon}$,
\quad\= by $(\c{\kon}\c{\kon})$ and $\mbox{($\c{\kon}$ {\it
nat})},$
\\[1ex]
\> \>
$=\;\b{\kon}{\rts}\!\cirk\!\c{\kon}\!\cirk\!\b{\kon}{\rts}\!\cirk\!\Dk\!\cirk\!\c{\kon}$,
\> by $\mbox{($\b{\kon}{}\Dk$)},$
\\[1ex]
\> \>
$=(\mj\:\kon\c{\kon})\cirk\!\b{\kon}{\rts}\!\cirk\!\Dk$,\quad with
$\mbox{($\Dk$ {\it nat})}$ and $(\b{\kon}{}\c{\kon}),$
\\[1ex]
\> \> $=\mj\:\kon\Sk$, \quad by $\mbox{($\b{\kon}{}\Dk$)}$,
\end{tabbing}

\noindent and for the second equation we have

\begin{tabbing}
\mbox{\hspace{2em}}\= $\b{\kon}{\str}\!\cirk\!\Sk\;$\=
$=\;\b{\kon}{\str}\!\cirk\!\c{\kon}\!\cirk\!\b{\kon}{\str}\!\cirk(\mj\:\kon\Dk)$,
\quad with $\mbox{($\b{\kon}{}\Dk$)}$,
\\[1ex]
\> \> $=(\c{\kon}\kon\:\mj)\cirk\!\b{\kon}{\str}\!
\cirk(\mj\:\kon\c{\kon})\cirk(\mj\:\kon\Dk)$,\quad by
$(\b{\kon}{}\c{\kon})$,
\\[1ex]
\> \> $=\;\Sk\kon\:\mj$,\quad with $\mbox{($\b{\kon}{}\Dk\Sk$)}$.
\end{tabbing}

\noindent We derive analogously with the help of
$\mbox{($\b{\dis}{}\Sd$)}$ the equations

\begin{tabbing}
\mbox{\hspace{2em}}\= $\mbox{($\Dkp$ {\it nat})}$\quad\quad\=
$(f\kon\mj_{B\dis\neg B})\cirk\!\Dkp_{B,A}\;=\;\Dkp_{B,D}\!\cirk
f$.\kill

\mbox{\hspace{2em}}\> $\mbox{($\b{\dis}{}\Dd\Sd$)}$\>
$(\Dd_{B,A}\dis\:\mj_C)\cirk\!\b{\dis}{\str}_{A,B\kon\neg
B,C}\;=\mj_A\:\dis\Sd_{B,C}$,\index{bdisdeltadissigmadis@(${\b{\dis}{}\Dd\Sd}$)
equation}
\\[1ex]
\> $\mbox{($\b{\dis}{}\Dd$)}$\> $\Dd_{C,A\dis
B}\!\cirk\!\b{\dis}{\str}_{A,B,C\kon\neg
C}\;=\mj_A\:\dis\Dd_{C,B}$.\index{bdisdeltadis@(${\b{\dis}{}\Dd}$)
equation}
\end{tabbing}

The arrows ${\Dk_{B,A}:A\vdash A\kon(\neg B\dis B)}$ and
$\Sk_{B,A}:A\vdash(\neg B\dis B)\kon A$ are analogous to the
arrows of types ${A\vdash A\kon\top}$ and ${A\vdash\top\kon A}$
that one finds in monoidal categories. However, $\Dk_{B,A}$ and
$\Sk_{B,A}$ do not have inverses in \PNN. The equations
$\mbox{($\b{\kon}{}\Dk$)}$, $\mbox{($\b{\kon}{}\Dk\Sk$)}$,
$\mbox{($\b{\kon}{}\Sk$)}$ are analogous to equations that hold in
monoidal categories (see \cite{ML71}, Section VII.1, \cite{DP04},
Section 4.6). An analogous remark can be made for $\Sd_{B,A}$ and
$\Dd_{B,A}$.

We can also derive for \PNN\ the following equations by using
essentially $\mbox{($d\!\Sk$)}$ and $\mbox{($d\!\Dd$)}$:

\begin{tabbing}
\mbox{\hspace{2em}}\= $\mbox{($\Dkp$ {\it nat})}$\quad\quad\=
$(f\kon\mj_{B\dis\neg B})\cirk\!\Dkp_{B,A}\;=\;\Dkp_{B,D}\!\cirk
f$.\kill

\mbox{\hspace{2em}}\> $\mbox{($d^R\!\Dk$)}$\> $d^R_{C,B,\neg A\dis
A}\cirk\!\Dk_{A,C\dis B}\;$\=
$=\mj_C\:\dis\Dk_{A,B}$,\index{dassocrdeltacon@(${d^R\Dk}$)
equation}
\\[1ex]
\> $\mbox{($d^R\!\Sd$)}$\> $\Sd_{A,B\kon C}\!\cirk d^R_{A\kon\neg
A,B,C}$\>
$=\;\Sd_{A,B}\kon\:\mj_C$.\index{dassocrsigmadis@(${d^R\Sd}$)
equation}
\end{tabbing}

\noindent These two equations could replace $\mbox{($d\!\Sk$)}$
and $\mbox{($d\!\Dd$)}$ for defining \PNN. The analogues of the
equations $\mbox{($d\!\Sk$)}$, $\mbox{($d\!\Dd$)}$,
$\mbox{($d^R\!\Dk$)}$ and $\mbox{($d^R\!\Sd$)}$ may be found in
\cite{CS97} (Section 2.1), where they are assumed for linearly
(alias weakly) distributive categories with negation (cf.\
\cite{DP04}, Section 7.9).

It is easy to infer that in \PNN\ we have analogues of the
equations $\mbox{($\b{\kon}{}\Dk$)}$,
$\mbox{($\b{\kon}{}\Dk\Sk$)}$, $\mbox{($\b{\kon}{}\Sk$)}$,
$\mbox{($\b{\dis}{}\Sd$)}$, $\mbox{($\b{\dis}{}\Dd\Sd$)}$,
$\mbox{($\b{\dis}{}\Dd$)}$, $\mbox{($d\!\Sk$)}$,
$\mbox{($d\!\Dd$)}$, $\mbox{($d^R\!\Dk$)}$ and
$\mbox{($d^R\!\Sd$)}$ obtained by replacing $\Xx$ by $\Xxp\!$, and
the indices of the form ${\neg B\dis B}$ and ${B\kon\neg B}$ by
${B\dis\neg B}$ and ${\neg B\kon B}$ respectively. For example, we
have

\begin{tabbing}
\mbox{\hspace{2em}}\= $\mbox{($\Dkp$ {\it nat})}$\quad\quad\=
$(f\kon\mj_{B\dis\neg B})\cirk\!\Dkp_{B,A}\;=\;\Dkp_{B,D}\!\cirk
f$.\kill

\> $\mbox{($\b{\kon}{}\Dkp$)}$\> $\b{\kon}{\rts}_{A,B,C\dis\neg
C}\!\cirk\!\Dkp_{C,A\kon
B}\;=\mj_A\:\kon\Dkp_{C,B}$.\index{bcondeltaconprime@(${\b{\kon}{}\Dkp}$)
equation}
\end{tabbing}

We can also derive for \PNN\ the following equations by using
essentially $\mbox{($\Sd\Dk$)}$ and $\mbox{($\Sdp\Dkp$)}$:

\begin{tabbing}
\mbox{\hspace{2em}}\= $\mbox{($\Dkp$ {\it nat})}$\quad\quad\=
$(f\kon\mj_{B\dis\neg B})\cirk\!\Dkp_{B,A}\;=\;\Dkp_{B,D}\!\cirk
f$.\kill

\> $\mbox{($\Ddp\Skp$)}$\> $\Ddp_{A,A}\!\cirk d^R_{A,\neg
A,A}\cirk\!\Skp_{A,A}\;=\mj_A$,\index{Deltadisprimesigmaconprime@(${\Ddp\Skp}$)
equation}
\\[1ex]
\> $\mbox{($\Dd\Sk$)}$\> $\Dd_{A,\neg A}\!\cirk d^R_{\neg A,A,\neg
A}\cirk\!\Sk_{A,\neg A}\;=\mj_{\neg
A}$.\index{Deltadissigmacon@(${\Dd\Sk}$) equation}
\end{tabbing}

\noindent These two equations could replace $\mbox{($\Sd\Dk$)}$
and $\mbox{($\Sdp\Dkp$)}$ for defining \PNN. The equations
$\mbox{($\Sd\Dk$)}$, $\mbox{($\Sdp\Dkp$)}$, $\mbox{($\Ddp\Skp$)}$
and $\mbox{($\Dd\Sk$)}$ are related to the triangular equations of
an adjunction (see \cite{ML71}, Section IV.1; see also the next
section). The analogues of these equations may be found in
\cite{CS97} (Section~4).

A {\it proof-net} category\index{proof-net category} is a category
with two biendofunctors $\kon$ and $\dis$, a unary operation
$\neg$ on objects, and the natural transformations
$\b{\kon}{\str}$, $\b{\kon}{\rts}$, $\b{\dis}{\str}$,
$\b{\dis}{\rts}$, $\c{\kon}$, $\c{\dis}$, $d$, $\Dk$ and $\Sd$
that satisfy the equations $(\b{\xi}{}\!5)$,
$(\b{\xi}{}\b{\xi}{})$, $\ldots$ , $\mbox{($\Sdp\Dkp$)}$ of \PNN.
The category \PNN\ is up to isomorphism the free proof-net
category generated by the set of letters \Pe\ (the set \Pe\ may be
understood as a discrete category).

If $\beta$ is a primitive arrow term of \PNN\ except $\mj_B$, then
we call $\beta$-{\it terms}\index{beta term@$\beta$-term} of \PNN\
the set of arrow terms defined inductively as follows: $\beta$ is
a $\beta$-term; if $f$ is a $\beta$-term, then for every $A$ in
${\eL_{\kon,\dis}}$ we have that ${\mj_A\ks f}$ and ${f\ks
\mj_A}$, where ${\!\ks\!\in\{\kon,\dis\}}$, are $\beta$-terms.

In a $\beta$-term the subterm $\beta$ is called the {\it
head}\index{head} of this $\beta$-term. For example, the head of
the $\b{\kon}{\str}_{B,C,D}$-term
${\mj_A\kon(\b{\kon}{\str}_{B,C,D}\dis\:\mj_E)}$ is
$\b{\kon}{\str}_{B,C,D}$.

We define $\mj$-\emph{terms}\index{identity term@identity term
$\mj$-term} as $\beta$-terms by replacing $\beta$ in the
definition above by $\mj_B$. So $\mj$-terms are headless.

An arrow term of the form ${f_n \cirk \ldots \cirk f_1}$, where
$n\geq 1$, with parentheses tied to $\cirk$ associated
arbitrarily, such that for every ${i\in\{1,\ldots,n\}}$ we have
that $f_i$ is composition-free is called {\it
factorized}.\index{factorized arrow term} In a factorized arrow
term ${f_n\cirk\ldots\cirk f_1}$ the arrow terms $f_i$ are called
{\it factors}.\index{factor} A factor that is a $\beta$-term for
some $\beta$ is called a {\it headed}\index{headed factor} factor.
A factorized arrow term is called {\it headed}\index{headed
factorized arrow term} when each of its factors is either headed
or a $\mj$-term. A factorized arrow term ${f_n\cirk\ldots\cirk
f_1}$ is called {\it developed}\index{developed factorized arrow
term} when $f_1$ is a $\mj$-term and if $n>1$, then every factor
of ${f_n\cirk\ldots\cirk f_2}$ is headed. It is sometimes useful
to write the factors of a headed arrow term one above the other,
as it is done for example in Figure 1 at the end of \S 6.

By using the categorial equations $({\mbox{{\it cat}~1}})$ and
$({\mbox{{\it cat}~2}})$ and bifunctorial equations we can easily
prove by induction on the length of $f$ the following lemma.

\prop{Development Lemma}{For every arrow term $f$ there is a
developed arrow term $f'$ such that $f=f'$ in
\PNN.}\index{Development Lemma}

\noindent Analogous definitions of $\beta$-term and developed
arrow term can be given for \mds, and an analogous Development
Lemma can be proved for \mds.

\section{\large\bf The category \emph{Br}}

We are now going to introduce a category called \emph{Br}, which
will serve to prove our main coherence result for proof-net
categories. We will show that there is a faithful functor from
\PNN\ to \emph{Br}. The name of the category \emph{Br} comes from
``Brauerian''. The arrows of this category correspond to graphs,
or diagrams, that were introduced in \cite{B37} in connection with
Brauer algebras. Analogous graphs were investigated in
\cite{EK66}, and in \cite{KML71} Kelly\index{Kelly@Kelly, G.M.}
and Mac Lane\index{Mac Lane@Mac Lane, S.} relied on them to prove
their coherence result for symmetric monoidal closed categories.

Let \eM\ be a set whose subsets are denoted by $X$, $Y$,
$Z$,~$\ldots$ For ${i\in\{s,t\}}$ (where $s$ stands for ``source''
and $t$ for ``target''),\index{s source@$s$ source}\index{t
target@$t$ target} let $\eM^i$ be a set in one-to-one
correspondence with \eM, and let ${i\!:\eM\str \eM^i}$ be a
bijection. Let $X^i$ be the subset of $\eM^i$ that is the image of
the subset $X$ of \eM\ under $i$. If ${u\in\eM}$, then we use
$u_i$ as an abbreviation for $i(u)$. We assume also that \eM,
$\eM^s$ and $\eM^t$ are mutually disjoint.

For ${X,Y\subseteq \eM}$, let a {\it split relation}\index{split
relation} of \eM\ be a triple ${\langle R,X,Y\rangle}$ such that
${R\subseteq(X^s\cup Y^t)^2}$. The set ${X^s\cup Y^t}$ may be
conceived as the disjoint union of $X$ and $Y$. We denote a split
relation ${\langle R,X,Y\rangle}$ more suggestively by
${R\!:X\vdash Y}$.

A split relation ${R\!:X\vdash Y}$ is a {\it split
equivalence}\index{split equivalence} when $R$ is an equivalence
relation. We denote by part($R$) the partition of ${X_s\cup Y_t}$
corresponding to the split equivalence ${R\!:X\vdash Y}$.

A split equivalence ${R\!:X\vdash Y}$ is {\it
Brauerian}\index{Brauerian split equivalence} when every member of
part($R$) is a two-element set. For ${R\!:X\vdash Y}$ a Brauerian
split equivalence, every member of part($R$) is either of the form
${\{u_s,v_t\}}$, in which case it is called a {\it
transversal},\index{transversal} or of the form ${\{u_s,v_s\}}$,
in which case it is called a {\it cup},\index{cup} or, finally, of
the form ${\{u_t,v_t\}}$, in which case it is called a {\it
cap}.\index{cocacap@cap}

For ${X,Y,Z\in\eM}$, we want to define the composition ${P\ast
R\!:X\vdash Z}$ of the split relations ${R\!:X\vdash Y}$ and
${P\!:Y\vdash Z}$ of \eM. For that we need some auxiliary notions.

For ${X,Y\subseteq \eM}$, let the function ${\varphi^s\!: X\cup
Y^t\str X^s\cup Y^t}$ be defined by

\[
\varphi^s(u)=\left\{
\begin{array}{ll}
u_s & \mbox{\rm{if }} u\in X
\\
u       & \mbox{\rm{if }} u\in Y^t,
\end{array}
\right.
\]

\noindent and let the function ${\varphi^t\!:X^s\cup Y\str X^s\cup
Y^t}$ be defined by

\[
\varphi^t(u)=\left\{
\begin{array}{ll}
u            & \mbox{\rm{if }} u\in X^s
\\
u_t       & \mbox{\rm{if }} u\in Y.
\end{array}
\right.
\]

For a split relation ${R\!:X\vdash Y}$, let the two relations
${R^{-s}\subseteq(X\cup Y^t)^2}$ and ${R^{-t}\subseteq(X^s\cup
Y)^2}$ be defined by

\[
(u,v)\in R^{-i}\quad\mbox{\rm{iff}}\quad
(\varphi^i(u),\varphi^i(v))\in R
\]

\noindent for ${i\in\{s,t\}}$. Finally, for an arbitrary binary
relation $R$, let $\mbox{\rm Tr}(R)$ be the transitive closure of
$R$.

Then we define ${P\ast R}$ by

\[
P\ast R=_{\df}\mbox{\rm Tr}(R^{-t}\cup P^{-s})\cap(X^s\cup Z^t)^2.
\]

\noindent It is easy to conclude that ${P\ast R\!:X\vdash Z}$ is a
split relation of \eM, and that if ${R\!:X\vdash Y}$ and
${P\!:Y\vdash Z}$ are (Brauerian) split equivalences, then ${P\ast
R}$ is a (Brauerian) split equivalence.

We now define the category \emph{Br}.\index{Br@\emph{Br} category}
The objects of \emph{Br} are the members of the set of finite
ordinals {\boldmath$N$}.\index{N@{\boldmath$N$ finite
ordinals}}\index{finite ordinals} (We have ${0=\emptyset}$ and
${n\pl 1=n\cup\{n\}}$, while {\boldmath$N$} is the ordinal
$\omega$.) The arrows of \emph{Br} are the Brauerian split
equivalences ${R\!:m\vdash n}$ of {\boldmath$N$}. The identity
arrow ${\mj_n\!:n\vdash n}$ of \emph{Br}\index{identity arrow of
Br@identity arrow of \emph{Br} $\mj_n$} is the Brauerian split
equivalence such that

\[
\mbox{\rm{part}}(\mj_n)=\{\{m_s,m_t\}\mid m<n\}.
\]

\noindent Composition in \emph{Br} is the operation $\ast$ defined
above.

That \emph{Br} is indeed a category (i.e.\ that $\ast$ is
associative and that $\mj_n$ is an identity arrow) is proved in
\cite{DP03a} and \cite{DP03b}. This proof is obtained via an
isomorphic representation of \emph{Br} in the category
\emph{Rel},\index{Rel@\emph{Rel} category} whose objects are the
finite ordinals and whose arrows are all the relations between
these objects. Composition in \emph{Rel} is the ordinary
composition of relations. A direct formal proof would be more
involved, though what we have to prove is rather clear if we
represent Brauerian split equivalences geometrically (as this is
done in \cite{B37}, \cite{EK66}, and also in categories of
tangles; see \cite{K95}, Chapter 12, and references therein).

For example, for ${R\subseteq(3^s\cup 9^t)^2}$ and
${P\subseteq(9^s\cup 1^t)^2}$ such that

\begin{tabbing}
\hspace{2em}\=$\mbox{\rm{part}}(R)\;$\=$=\{\{0_s,0_t\},\{1_s,3_t\},\{2_s,6_t\}\}\cup
\{\{n_t,(n\pl 1)_t\}\mid n\in\{1,4,7\}\}$,
\\*[1.5ex]
\>$\mbox{\rm{part}}(P)$\>$=\{\{2_s,0_t\}\}\cup\{\{n_s,(n\pl
1)_s\}\mid n\in\{0,3,5,7\}\}$,
\end{tabbing}

\noindent the composition ${P\ast R\subseteq(3^s\cup 1^t)^2}$, for
which we have

\[
\mbox{\rm{part}}(P\ast R)={\{\{0_s,0_t\}, \{1_s,2_s\}\}},
\]

\noindent is obtained from the following diagram:

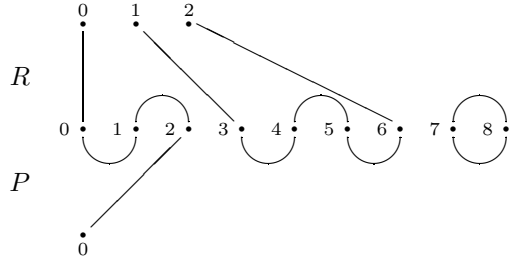
\begin{center}
\begin{picture}(160,115)

\put(3,23){\line(1,1){34}} \put(0,63){\line(0,1){34}}
\put(57,63){\line(-1,1){34}} \put(117,63){\line(-2,1){74}}

\put(0,20){\circle*{2}} \put(0,60){\circle*{2}}
\put(20,60){\circle*{2}} \put(40,60){\circle*{2}}
\put(60,60){\circle*{2}} \put(80,60){\circle*{2}}
\put(100,60){\circle*{2}} \put(120,60){\circle*{2}}
\put(140,60){\circle*{2}} \put(160,60){\circle*{2}}
\put(0,100){\circle*{2}} \put(20,100){\circle*{2}}
\put(40,100){\circle*{2}}

\put(10,57){\oval(20,20)[b]} \put(30,63){\oval(20,20)[t]}
\put(70,57){\oval(20,20)[b]} \put(110,57){\oval(20,20)[b]}
\put(90,63){\oval(20,20)[t]} \put(150,57){\oval(20,20)[b]}
\put(150,63){\oval(20,20)[t]}

\put(0,17){\makebox(0,0)[t]{\scriptsize$0$}}
\put(-5,60){\makebox(0,0)[r]{\scriptsize$0$}}
\put(15,60){\makebox(0,0)[r]{\scriptsize$1$}}
\put(35,60){\makebox(0,0)[r]{\scriptsize$2$}}
\put(55,60){\makebox(0,0)[r]{\scriptsize$3$}}
\put(75,60){\makebox(0,0)[r]{\scriptsize$4$}}
\put(95,60){\makebox(0,0)[r]{\scriptsize$5$}}
\put(115,60){\makebox(0,0)[r]{\scriptsize$6$}}
\put(135,60){\makebox(0,0)[r]{\scriptsize$7$}}
\put(155,60){\makebox(0,0)[r]{\scriptsize$8$}}
\put(0,103){\makebox(0,0)[b]{\scriptsize$0$}}
\put(20,103){\makebox(0,0)[b]{\scriptsize$1$}}
\put(40,103){\makebox(0,0)[b]{\scriptsize$2$}}

\put(-20,80){\makebox(0,0)[r]{$R$}}
\put(-20,40){\makebox(0,0)[r]{$P$}}

\end{picture}
\end{center}

\vspace{-2ex}

Every bijection $f$ from $X^s$ to $Y^t$ corresponds to a Brauerian
split equivalence ${R\!:X\vdash Y}$ such that the members of
part($R$) are of the form ${\{u,f(u)\}}$. The composition of such
Brauerian split equivalences, which correspond to bijections, is
then a simple matter: it amounts to composition of these
bijections. If in \emph{Br} we keep as arrows only such Brauerian
split equivalences, then we obtain a subcategory of \emph{Br}
isomorphic to the category \emph{Bij}\index{Bij@\emph{Bij}
category} whose objects are again the finite ordinals and whose
arrows are the bijections between these objects. The category
\emph{Bij} is a subcategory of the category
\emph{Rel}\index{Rel@\emph{Rel} category} (which played an
important role in \cite{DP04}), whose objects are the finite
ordinals and whose arrows are all the relations between these
objects. Composition in \emph{Bij} and \emph{Rel} is the ordinary
composition of relations. The category \emph{Rel} (which played an
important role in \cite{DP04}) is isomorphic to a subcategory of
the category whose arrows are split relations of finite ordinals,
of whom \emph{Br} is also a subcategory.

We define a functor $G$\index{G functor@$G$ functor} from \PNN\ to
\emph{Br} in the following way. On objects, we stipulate that $GA$
is the number of occurrences of letters in $A$. (If $A$ has
${n=\{0,1,\ldots,n\!-\! 1\}}$ occurrences of letters, then the
first occurrence corresponds to $0$, the second to $1$, etc.) On
arrows, we have first that ${G\alpha}$ is an identity arrow of
\emph{Br} for $\alpha$ being $\mj_A$, $\b{\xi}{\str}_{A,B,C}$,
$\b{\xi}{\rts}_{A,B,C}$ and $d_{A,B,C}$, where
${\!\ks\!\in\{\kon,\dis\}}$.

Next, for ${i,j\in\{s,t\}}$, we have that ${\{m_i,n_j\}}$ belongs
to part(${G\!\c{\kon}_{A,B}}$) iff ${\{n_i,m_j\}}$ belongs to
part(${G\!\c{\dis}_{A,B}}$), iff $i$ is $s$ and $j$ is $t$, while
${m,n<GA\pl GB}$ and

\vspace{-1ex}

\[
(m\mn n\mn GA)(m\mn n\pl GB)=0.
\]

\noindent In the following example, we have ${G(p\dis
q)=2=\{0,1\}}$ and $G((q\dis\neg r)\dis q)$$=3=\{0,1,2\}$, and we
have the diagrams

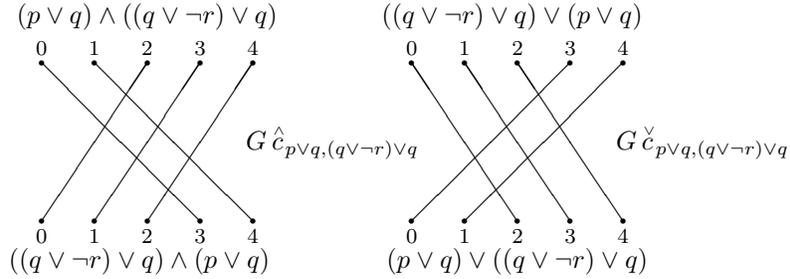
\vspace{2ex}

\begin{center}
\begin{picture}(300,95)
\put(10,10){\line(2,3){40}} \put(30,10){\line(2,3){40}}
\put(50,10){\line(2,3){40}} \put(70,10){\line(-1,1){60}}
\put(90,10){\line(-1,1){60}}

\put(150,10){\line(1,1){60}} \put(170,10){\line(1,1){60}}
\put(190,10){\line(-2,3){40}} \put(210,10){\line(-2,3){40}}
\put(230,10){\line(-2,3){40}}

\put(10,10){\circle*{2}} \put(30,10){\circle*{2}}
\put(50,10){\circle*{2}} \put(70,10){\circle*{2}}
\put(90,10){\circle*{2}} \put(150,10){\circle*{2}}
\put(170,10){\circle*{2}} \put(190,10){\circle*{2}}
\put(210,10){\circle*{2}} \put(230,10){\circle*{2}}

\put(10,70){\circle*{2}} \put(30,70){\circle*{2}}
\put(50,70){\circle*{2}} \put(70,70){\circle*{2}}
\put(90,70){\circle*{2}} \put(150,70){\circle*{2}}
\put(170,70){\circle*{2}} \put(190,70){\circle*{2}}
\put(210,70){\circle*{2}} \put(230,70){\circle*{2}}

\put(10,7){\makebox(0,0)[t]{\footnotesize 0}}
\put(30,7){\makebox(0,0)[t]{\footnotesize 1}}
\put(50,7){\makebox(0,0)[t]{\footnotesize 2}}
\put(70,7){\makebox(0,0)[t]{\footnotesize 3}}
\put(90,7){\makebox(0,0)[t]{\footnotesize 4}}
\put(150,7){\makebox(0,0)[t]{\footnotesize 0}}
\put(170,7){\makebox(0,0)[t]{\footnotesize 1}}
\put(190,7){\makebox(0,0)[t]{\footnotesize 2}}
\put(210,7){\makebox(0,0)[t]{\footnotesize 3}}
\put(230,7){\makebox(0,0)[t]{\footnotesize 4}}

\put(10,73){\makebox(0,0)[b]{\footnotesize 0}}
\put(30,73){\makebox(0,0)[b]{\footnotesize 1}}
\put(50,73){\makebox(0,0)[b]{\footnotesize 2}}
\put(70,73){\makebox(0,0)[b]{\footnotesize 3}}
\put(90,73){\makebox(0,0)[b]{\footnotesize 4}}
\put(150,73){\makebox(0,0)[b]{\footnotesize 0}}
\put(170,73){\makebox(0,0)[b]{\footnotesize 1}}
\put(190,73){\makebox(0,0)[b]{\footnotesize 2}}
\put(210,73){\makebox(0,0)[b]{\footnotesize 3}}
\put(230,73){\makebox(0,0)[b]{\footnotesize 4}}

\put(120,40){\makebox(0,0){$G\!\c{\kon}_{p\dis q,(q\dis\neg r)\dis
q}$}}

\put(260,40){\makebox(0,0){$G\!\c{\dis}_{p\dis q,(q\dis\neg r)\dis
q}$}}

\put(50,85){\makebox(0,0)[b]{$(p\dis q)\kon((q\dis\neg r)\dis
q)$}}

\put(47,0){\makebox(0,0)[t]{$((q\dis\neg r)\dis q)\kon(p\dis q)$}}

\put(188,85){\makebox(0,0)[b]{$((q\dis\neg r)\dis q)\dis(p\dis
q)$}}

\put(190,0){\makebox(0,0)[t]{$(p\dis q)\dis((q\dis\neg r)\dis
q)$}}

\end{picture}
\end{center}

\vspace{2ex}

We have that ${\{m_i,n_j\}}$ belongs to part(${G\!\Dk_{B,A}}$) iff
either

\nav{}{$i$ is $s$ and $j$ is $t$, while ${m,n<GA}$ and ${m=n}$,
or}

\vspace{-1ex}

\nav{}{$i$ and $j$ are both $t$, while ${m,n\in\{GA,\ldots,GA\pl
2GB\mn 1\}}$ and ${|m\mn n|=GB}$.}

\noindent In the following example, for $A$ being ${(q\dis\neg
r)\dis q}$ and $B$ being ${p\dis q}$, we have

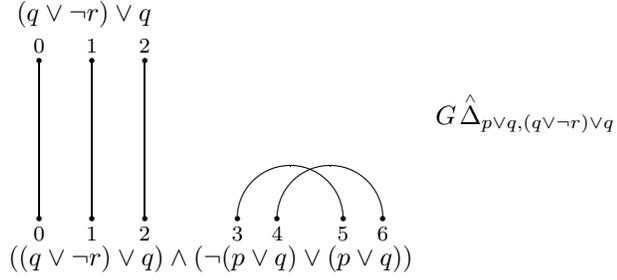
\vspace{2ex}

\begin{center}
\begin{picture}(180,107)
\put(10,20){\line(0,1){60}} \put(30,20){\line(0,1){60}}
\put(50,20){\line(0,1){60}}

\put(10,20){\circle*{2}} \put(30,20){\circle*{2}}
\put(50,20){\circle*{2}} \put(85,20){\circle*{2}}
\put(100,20){\circle*{2}} \put(125,20){\circle*{2}}
\put(140,20){\circle*{2}}

\put(10,80){\circle*{2}} \put(30,80){\circle*{2}}
\put(50,80){\circle*{2}}

\put(10,17){\makebox(0,0)[t]{\footnotesize 0}}
\put(30,17){\makebox(0,0)[t]{\footnotesize 1}}
\put(50,17){\makebox(0,0)[t]{\footnotesize 2}}
\put(85,17){\makebox(0,0)[t]{\footnotesize 3}}
\put(100,17){\makebox(0,0)[t]{\footnotesize 4}}
\put(125,17){\makebox(0,0)[t]{\footnotesize 5}}
\put(140,17){\makebox(0,0)[t]{\footnotesize 6}}

\put(10,83){\makebox(0,0)[b]{\footnotesize 0}}
\put(30,83){\makebox(0,0)[b]{\footnotesize 1}}
\put(50,83){\makebox(0,0)[b]{\footnotesize 2}}

\put(105,20){\oval(40,40)[t]} \put(120,20){\oval(40,40)[t]}

\put(160,60){\makebox(0,0)[l]{$G\!\Dk_{p\dis q,(q\dis\neg r)\dis
q}$}}

\put(75,10){\makebox(0,0)[t]{$((q\dis\neg r)\dis q)\kon(\neg(p\dis
q) \dis(p\dis q))$}}

\put(27,95){\makebox(0,0)[b]{$(q\dis\neg r)\dis q$}}

\end{picture}
\end{center}

\vspace{2ex}

We have that ${\{m_i,n_j\}}$ belongs to part(${G\!\Sd_{B,A}}$) iff
either

\nav{}{$i$ is $s$ and $j$ is $t$, while ${m\in\{2GB,\ldots,2GB\pl
GA\mn 1\}}$, ${n<GA}$ and ${m\mn 2GB=n}$, or}

\vspace{-1ex}

\nav{}{$i$ and $j$ are both $s$, while ${m,n<2GB}$ and ${|m\mn
n|=GB}$.}

\noindent For $A$ and $B$ being as in the previous example, we
have

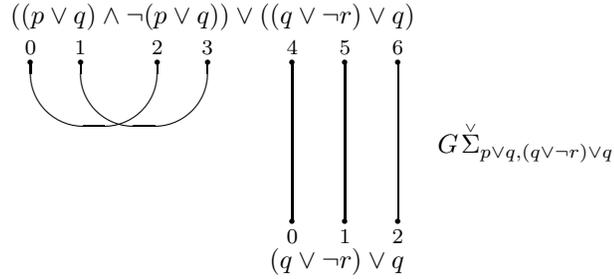
\begin{center}
\begin{picture}(180,110)
\put(100,20){\line(0,1){60}} \put(120,20){\line(0,1){60}}
\put(140,20){\line(0,1){60}}

\put(1,80){\circle*{2}} \put(20,80){\circle*{2}}
\put(49,80){\circle*{2}} \put(68,80){\circle*{2}}
\put(100,80){\circle*{2}} \put(120,80){\circle*{2}}
\put(140,80){\circle*{2}}

\put(100,20){\circle*{2}} \put(120,20){\circle*{2}}
\put(140,20){\circle*{2}}

\put(1,83){\makebox(0,0)[b]{\footnotesize 0}}
\put(20,83){\makebox(0,0)[b]{\footnotesize 1}}
\put(49,83){\makebox(0,0)[b]{\footnotesize 2}}
\put(68,83){\makebox(0,0)[b]{\footnotesize 3}}
\put(100,83){\makebox(0,0)[b]{\footnotesize 4}}
\put(120,83){\makebox(0,0)[b]{\footnotesize 5}}
\put(140,83){\makebox(0,0)[b]{\footnotesize 6}}

\put(100,17){\makebox(0,0)[t]{\footnotesize 0}}
\put(120,17){\makebox(0,0)[t]{\footnotesize 1}}
\put(140,17){\makebox(0,0)[t]{\footnotesize 2}}

\put(25,80){\oval(48,48)[b]} \put(44,80){\oval(48,48)[b]}

\put(155,50){\makebox(0,0)[l]{$G\!\Sd_{p\dis q,(q\dis\neg r)\dis
q}$}}

\put(117,10){\makebox(0,0)[t]{$(q\dis\neg r)\dis q$}}

\put(70,95){\makebox(0,0)[b]{$((p\dis q)\kon\neg(p\dis q))\dis
((q\dis\neg r)\dis q)$}}

\end{picture}
\end{center}

Let ${G(f\cirk g)=Gf\ast Gg}$. To define ${G(f\ks g)}$, for
${\!\ks\!\in\{\kon,\dis\}}$,  we need an auxiliary notion.

Suppose $b_X$ is a bijection from $X$ to $X_1$ and $b_Y$ a
bijection from $Y$ to $Y_1$. Then for ${R\subseteq(X^s\cup
Y^t)^2}$ we define ${R_{b_Y}^{b_X}\subseteq(X_1^s\cup Y_1^t)^2}$
by

\[(u_i,v_j)\in R_{b_Y}^{b_X}
\quad\mbox{\rm{iff}}\quad (i(b_U^{-1}(u)),j(b_V^{-1}(v)))\in R,
\]

\noindent where ${(i,U),(j,V)\in\{(s,X),(t,Y)\}}$.

If ${f\!:A\vdash D}$ and ${g\!:B\vdash E}$, then for
${\!\ks\!\in\{\kon,\dis\}}$ the set of ordered pairs ${G(f\ks g)}$
is

\[
Gf\cup Gg_{+GD}^{+GA}
\]

\noindent where ${+GA}$ is the bijection from $GB$ to $\{n\pl
GA\mid n\in GB\}$ that assigns ${n\pl GA}$ to $n$, and ${+GD}$ is
the bijection from $GE$ to $\{n\pl GD\mid n\in GE\}$ that assigns
${n\pl GD}$ to $n$.

It is not difficult to check that $G$ so defined is indeed a
functor from \PNN\ to \emph{Br}. For that, we determine by
induction on the length of derivation that for every equation
${f=g}$ of \PNN\ we have ${Gf=Gg}$ in \emph{Br}.

Consider, for example, the following diagram, which illustrates an
instance of $\mbox{($\Sd\Dk$)}$:

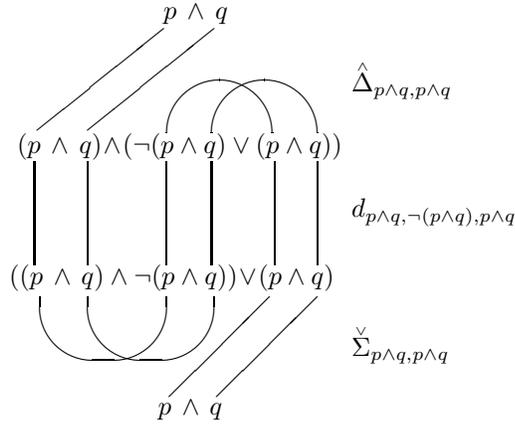
\begin{center}
\begin{picture}(150,180)

\put(51,20){\line(1,1){38}}

\put(69,20){\line(1,1){38}}

\put(0,70){\line(0,1){39}}

\put(20,70){\line(0,1){39}}

\put(50,70){\line(0,1){39}}

\put(67,70){\line(0,1){39}}

\put(91,70){\line(0,1){39}}

\put(107,70){\line(0,1){39}}

\put(1,121){\line(5,4){48}}

\put(20,121){\line(5,4){48}}

\put(26,58){\oval(48,48)[b]}

\put(44,58){\oval(48,48)[b]}

\put(70,120){\oval(40,40)[t]}

\put(87,120){\oval(40,40)[t]}

\put(120,40){\makebox(0,0)[l]{$\Sd_{p\kon q,p\kon q}$}}
\put(120,90){\makebox(0,0)[l]{$d_{p\kon q,\neg(p\kon q),p\kon
q}$}} \put(120,140){\makebox(0,0)[l]{$\Dk_{p\kon q,p\kon q}$}}

\put(59,15){\makebox(0,0){$p\,\kon\, q$}}

\put(52,65){\makebox(0,0){$((p\,\kon \,q)\kon\neg(p\kon
q))\!\dis\!(p\kon q)$}}

\put(55,115){\makebox(0,0){$(p\,\kon\, q)\!\kon\!(\neg(p\kon
q)\dis(p\kon q))$}}

\put(61,165){\makebox(0,0){$p\,\kon\, q$}}

\end{picture}
\end{center}

\vspace{-2ex}

\noindent This diagram shows that the equation
$\mbox{($\Sd\Dk$)}$, as well as the equation
$\mbox{($\Sdp\Dkp$)}$, which is illustrated by analogous diagrams,
is related to triangular equations of adjunctions (cf.\
\cite{D99}, Section 4.10). The triangular equations of adjunctions
are essentially about ``straightening a
serpentine'',\index{straightening a
serpentine}\index{smbserpentine@serpentine} and this straightening
is based on planar ambient isotopies of knot theory (cf.\
\cite{BZ85}, Section 1.A).

We have shown by this induction that \emph{Br} is a proof-net
category, and the existence of a structure-preserving functor $G$
from \PNN\ to \emph{Br} follows from the freedom of \PNN.

We can define analogously to $G$ a functor, which we also call
$G$,\index{G functor@$G$ functor} from the category \mds\ to
\emph{Br}. We just omit from the definition of $G$ above the
clauses involving $\Dk_{B,A}$ and $\Sd_{B,A}$. The image of \mds\
by $G$ in \emph{Br} is the subcategory of \emph{Br} isomorphic to
\emph{Bij}, which we mentioned above. The following is proved in
\cite{DP04} (Section 7.6).

\prop{\mds\ Coherence}{The functor $G$ from \mds\ to Br is
faithful.}\index{DS Coherence@\mds\ Coherence}

\noindent It follows immediately from this coherence result that
\mds\ is isomorphic to a subcategory of \PNN\ (cf.\ \cite{DP04},
Section 14.4).

Up to the end of \S 8 we will be occupied with proving the
following.

\prop{\PNN\ Coherence}{The functor $G$ from \PNN\ to Br is
faithful.}\index{PNN Coherence@\PNN\ Coherence}

\noindent For this proof, we must deal first with some preliminary
matters.

\section{\large\bf Some properties of \mds}

In this section we will prove some results about the category
\mds, which we will be use to ascertain that particular equations
hold in \PNN. We need these results also for the proof of \PNN\
Coherence.

First we introduce a definition. Suppose $x$ is the  \mbox{$n$-th}
occurrence of a letter (counting from the left) in a formula $A$
of ${\eL_{\neg,\kon,\dis}}$, and $y$ is the \mbox{$m$-th}
occurrence of the same letter in a formula $B$ of
${\eL_{\neg,\kon,\dis}}$. Then we say that $x$ and $y$ are {\it
linked} in an arrow ${f\!:A\vdash B}$ of \PNN\ when in the
partition part($Gf$) we have ${\{(n\mn 1)_s,(m\mn 1)_t\}}$ as a
member. (Note that to find the \mbox{$n$-th} occurrence, we count
starting from 1, but the ordinal ${n>0}$ is ${\{0,\ldots,n\mn
1\}}$.) We have an analogous definition of linked occurrences of
the same letter for \mds: we just replace ${\eL_{\neg,\kon,\dis}}$
by ${\eL_{\kon,\dis}}$ and \PNN\ by \mds.

It is easy to established by induction on the complexity of $f$
that for every arrow term ${f\!:A\vdash B}$ of \mds\ we have
${GA=GB}$. Moreover, every occurrence of letter in $A$ is linked
to exactly one occurrence of the same letter in $B$, and vice
versa. This is related to the fact that every arrow term
${f\!:A\vdash B}$ of \mds\ may be obtained by substituting letters
for letters out of an arrow term ${f'\!:A'\vdash B'}$ of \mds\
such that every letter occurs in $A'$ at most once, and the same
for $B'$ (see \cite{DP04}, Sections 3.3 and 7.6).

Suppose for Lemmata 1D and 2D below that ${f\!:A\vdash B}$ is an
arrow term of \mds\ such that $A$ has a subformula $D$ in which
$\kon$ does not occur and $B$ has a subformula $D'$ in which
$\kon$ does not occur, and suppose that every occurrence of a
letter in $D$ is linked to an occurrence of a letter in $D'$ and
vice versa. Then we can prove the following.

\prop{Lemma 1D}{The source $A$ of $f$ is $D$ iff the target $B$ of
$f$ is $D'$.}

\noindent This follows from the fact, noted above, that ${GA=GB}$.
The arrow term $f$ in this case can have as subterms that are
primitive arrow terms only arrow terms of the forms $\mj_E$,
$\b{\dis}{\str}_{E,F,G}$, $\b{\dis}{\rts}_{E,F,G}$ or
$\c{\dis}_{E,F}$. We also have the following.

\prop{Lemma 2D}{If ${D\kon A'}$ or ${A'\kon D}$ is a subformula of
$A$, then ${D'\kon B'}$ or ${B'\kon D'}$ is a subformula of $B$
for some $B'$.}

\noindent This is easily proved by induction on the complexity of
the arrow term $f$, with the help of Lemma~1D.

Suppose for Lemmata 1C and 2C below that ${f\!:A\vdash B}$ is an
arrow term of \mds\ such that $B$ has a subformula $C$ in which
$\dis$ does not occur and $A$ has a subformula $C'$ in which
$\dis$ does not occur, and suppose that every occurrence of a
letter in $C$ is linked to an occurrence of a letter in $C'$ and
vice versa. Then we can prove the following duals of Lemmata 1D
and 2D, in an analogous manner.

\prop{Lemma 1C}{The target $B$ of $f$ is $C$ iff the source $A$ of
$f$ is $C'$.}

\vspace{-2ex}

\prop{Lemma 2C}{If ${C\dis B'}$ or ${B'\dis C}$ is a subformula of
$B$, then ${C'\dis A'}$ or ${A'\dis C'}$ is a subformula of $A$
for some $A'$.}

Suppose for the following lemma, which is a corollary of either
Lemma~2D or Lemma~2C, that ${f\!:A\vdash B}$ is an arrow term of
\mds\ such that an occurrence $x$ of a letter $p$ in $A$ is linked
to an occurrence $y$ of $p$ in $B$.

\prop{Lemma 2}{It is impossible that $A$ has a subformula ${x\kon
A'}$ or ${A'\kon x}$ and $B$ has a subformula ${y\dis B'}$ or
${B'\dis y}$.}

Suppose for Lemmata 3D, 3C, 3 and 4 below that ${f\!:A\vdash B}$
is an arrow term of \mds, and for ${i\in\{1,2\}}$ let $x_i$ in $A$
and $y_i$ in $B$ be occurrences of the letter $p_i$ linked in $f$
(here $p_1$ and $p_2$ may also be the same letter).

\prop{Lemma 3D}{If in $A$ we have a subformula ${A_1\dis A_2}$
such that $x_i$ occurs in $A_i$, then in $B$ we have a subformula
${B_1\dis B_2}$ or ${B_2\dis B_1}$ such that $y_i$ occurs in
$B_i$.}

\noindent This is easily proved by induction on the complexity of
the arrow term $f$. We prove analogously the following.

\prop{Lemma 3C}{If in $B$ we have a subformula ${B_1\kon B_2}$
such that $y_i$ occurs in $B_i$, then in $A$ we have a subformula
${A_1\kon A_2}$ or ${A_2\kon A_1}$ such that $x_i$ occurs in
$A_i$.}

As a corollary of either Lemma~3D or Lemma~3C we have the
following.

\prop{Lemma 3}{It is impossible that $A$ has a subformula
${x_1\dis x_2}$ or ${x_2\dis x_1}$ and $B$ has a subformula
${y_1\kon y_2}$ or ${y_2\kon y_1}$.}

The following lemma, dual to Lemma~3, is a corollary of Lemma~2.

\prop{Lemma 4}{It is impossible that $A$ has a subformula
${x_1\kon x_2}$ or ${x_2\kon x_1}$ and $B$ has a subformula
${y_1\dis y_2}$ or ${y_2\dis y_1}$.}

\noindent Lemma~3 is related to the acyclicity condition of proof
nets, while Lemma~4 is related to the connectedness condition (see
\cite{DR89}).

Next we can prove the following lemma.

\prop{\emph{p-q-r} Lemma}{Let\index{pqr Lemma@\emph{p-q-r} Lemma}
${f\!:A\vdash B}$ be an arrow of \mds, let $x_i$ for
${i\in\{1,2,3\}}$ be occurrences of the letters $p$, $q$ and $r$,
respectively, in $A$, and let $y_i$ be occurrences of the letters
$p$, $q$ and $r$, respectively, in $B$, such that $x_i$ and $y_i$
are linked in $f$. Let, moreover, ${x_2\dis x_3}$ be a subformula
of $A$ and ${y_1\kon y_2}$ a subformula of $B$. Then there is a
$d_{p,q,r}$-term ${h\!:A'\vdash B'}$ such that $x_i'$ are
occurrences of the letters $p$, $q$ and $r$, respectively, in the
source ${p\kon(q\dis r)}$ of the head of $h$ and $y_i'$ are
occurrences of the letters $p$, $q$ and $r$, respectively, in the
target ${(p\kon q)\dis r}$ of the head of $h$, such that for some
arrows ${f_x\!:A\vdash A'}$ and ${f_y\!:B'\vdash B}$ of \mds\ we
have ${f=f_y\cirk h\cirk f_x}$ in \mds, and $x_i$ is linked to
$x_i'$ in $f_x$, while $y_i'$ is linked to $y_i$ in $f_y$.}

\dkz The proof of this lemma, of which we give just a sketch,
relies on a cut-elimination and related results of \cite{DP04}
(Sections 7.7-8). We first find in the category \mgds\ introduced
in \cite{DP04} (Section 7.7) a cut-free Gentzen term
${f'\!:X\vdash Y}$, which corresponds to $f$, by the relationship
that exists between \mds\ and \mgds. According to the equations at
the beginning of Section 7.8 of \cite{DP04}, which are used for
the proof of the Invertibility Lemmata in the same section, in
\mgds\ we have the equation ${f'=f''}$ for a Gentzen term $f''$
that has as a subterm either
${\kon_{p,q}(\mj_p,\dis_{q,r}(\mj_q,\mj_r))}$ or
${\dis_{q,r}(\kon_{p,q}(\mj_p,\mj_q),\mj_r)}$ both of type
${p\kon(q\dis r)\vdash (p\kon q)\dis r}$. By the relationship that
exists between \mds\ and \mgds, we can find starting from $f''$ an
arrow term ${f_y\cirk h\cirk f_x}$ equal to $f$ in \mds, which
satisfies the conditions of the lemma. \qed

\vspace{2ex}

The full force of the Cut-Elimination Theorem of Section 7.7 of
\cite{DP04} is not essential for this proof, but applying this
theorem simplifies the proof.

\section{\large\bf The category \PN}

We now introduce a category called \PN,\index{PN@\PN\ category}
which is equivalent to \PNN. In the objects of \PN, the negation
connective $\neg$ will be prefixed only to letters, and hence
$\Dk_{B,A}$ and $\Sd_{B,A}$ will be primitive only for the crown
index $B$ being a letter. Here is the formal definition of \PN.

For \Pe\ being the set of letters that we used to generate
${\eL_{\kon,\dis}}$ and ${\eL_{\neg,\kon,\dis}}$ in \S\S 2-3, let
$\Pe^{\neg}$\index{Pa@$\Pe^{\neg}$ set} be the set ${\{\neg p\mid
p\in\Pe\}}$. The objects of \PN\ are the formulae of the
propositional language ${\eL_{\kon,\dis}^{\neg
p}}$\index{Lcondisneg@${\eL_{\kon,\dis}^{\neg p}}$ language}
generated from ${\Pe\cup\Pe^{\neg}}$ with the binary connectives
$\kon$ and $\dis$. To define the arrow terms of \PN, in the
inductive definition we had for the arrow terms of \mds\ we assume
in addition that for every formula $A$ of ${\eL_{\kon,\dis}^{\neg
p}}$ and every letter $p$

\begin{tabbing}
\centerline{$\Dk_{p,A}:A\vdash A\kon(\neg p\dis
p)$,}\index{Deltacon@$\Dk$ arrow}
\\*[1ex]
\centerline{$\Sd_{p,A}:(p\kon\neg p)\dis A\vdash
A$,}\index{Sigmadis@$\Sd$ arrow}
\end{tabbing}

\noindent are primitive arrow terms of \PN.

To define the arrows of \PN, we assume as additional equations in
the inductive definition we had for the equations of \mds\ the
\PN\ equations of \S 3 restricted to the arrow terms $\Dk_{p,A}$
and $\Sd_{p,A}$. This means that in $\mbox{($\Dk$ {\it nat})}$ and
$\mbox{($\Sd$ {\it nat})}$ the crown index $B$ will be $p$, in
$\mbox{($\b{\kon}{}\Dk$)}$ and $\mbox{($\b{\dis}{}\Sd$)}$ the
crown index $C$ will be $p$, and in $\mbox{($d\!\Sk$)}$,
$\mbox{($d\!\Dd$)}$, $\mbox{($\Sd\Dk$)}$ and $\mbox{($\Sdp\Dkp$)}$
the crown index $A$ will be $p$. We define $\Sk_{p,A}$,
$\Dd_{p,A}$, $\Dkp_{p,A}$, $\Sdp_{p,A}$, $\Skp_{p,A}$ and
$\Ddp_{p,A}$ for \PN\ as they were defined in \PNN\ in terms of
$\Dk_{p,A}$ and $\Sd_{p,A}$.

The following equations of \PN, and hence also of \PNN, which we
call {\it stem-increasing} equations,\index{stem-increasing
equations} enable us to have in developed arrow terms only
$\Dk_{A,B}$-terms and $\Sd_{A,B}$-terms that coincide with their
heads:

\begin{tabbing}
\mbox{\hspace{0em}}\= $\mbox{($\mj\:\kon\Dk$)}$\quad\=
$\mj_A\:\kon\Dk_{p,B}$ \= = \= $\b{\kon}{\rts}_{A,B,\neg p\dis
p}\!\cirk\!\Dk_{p,A\kon B}$, \quad by $\mbox{($\b{\kon}{}\Dk$)}$,
\\*[1.5ex]
\> $\mbox{($\Dk\kon\:\mj$)}$\> $\Dk_{p,B}\kon\:\mj_A$ \> = \=
$\c{\kon}_{A,B\kon(\neg p\dis p)}\!\cirk\!\b{\kon}{\rts}_{A,B,\neg
p\dis p} \!\cirk(\c{\kon}_{B,A}\kon\:\mj_{\neg p\dis
p})\cirk\!\Dk_{p,B\kon A}$,
\\*[1ex]
\` by $(\c{\kon}\c{\kon})$, $\mbox{($\c{\kon}$ {\it nat})}$,
$\mbox{($\mj\:\kon\Dk$)}$ and $\mbox{($\Dk$ {\it nat})}$,
\\[2ex]
\> $\mbox{($\mj\:\dis\Dk$)}$\> $\mj_A\:\dis\Dk_{p,B}$ \> = \>
$d^R_{A,B,\neg p\dis p}\cirk\!\Dk_{p,A\dis B}$, \quad by
$\mbox{($d^R\!\Dk$)}$,
\\*[1.5ex]
\> $\mbox{($\Dk\dis\:\mj$)}$\> $\Dk_{p,B}\dis\:\mj_A$ \> = \>
$\c{\dis}_{B\kon(\neg p\dis p),A}\!\cirk d^R_{A,B,\neg p\dis
p}\cirk (\c{\dis}_{A,B}\kon\:\mj_{\neg p\dis
p})\cirk\!\Dk_{p,B\dis A}$,
\\*[1ex]
\` by $(\c{\dis}\c{\dis})$, $\mbox{($\c{\dis}$ {\it nat})}$,
$\mbox{($\mj\:\dis\Dk$)}$ and $\mbox{($\Dk$ {\it nat})}$,
\\[2ex]
\> $\mbox{($\Sd\dis\:\mj$)}$\> $\Sd_{p,B}\dis\:\mj_A$ \> = \>
$\Sd_{p,B\dis A}\!\cirk\!\b{\dis}{\rts}_{p\kon\neg p,B,A}$, \quad
by $\mbox{($\b{\dis}{}\Sd$)}$,
\\*[1.5ex]
\> $\mbox{($\mj\:\dis\Sd$)}$\> $\mj_A\:\dis\Sd_{p,B}$\> = \>
$\Sd_{p,A\dis B}\!\cirk(\mj_{p\kon\neg
p}\dis\c{\dis}_{A,B})\cirk\! \b{\dis}{\rts}_{p\kon\neg
p,B,A}\!\cirk\!\c{\dis}_{(p\kon\neg p)\dis B,A}$,
\\*[1ex]
\` by $(\c{\dis}\c{\dis})$, $\mbox{($\c{\dis}$ {\it nat})}$,
$\mbox{($\Sd\dis\:\mj$)}$ and $\mbox{($\Sd$ {\it nat})}$,
\\[2ex]
\> $\mbox{($\Sd\kon\:\mj$)}$\> $\Sd_{p,B}\kon\:\mj_A$ \> = \>
$\Sd_{p,B\kon A}\!\cirk d^R_{p\kon\neg p,B,A}$, \quad by
$\mbox{($d^R\!\Sd$)}$,
\\[1.5ex]
\> $\mbox{($\mj\:\kon\Sd$)}$\> $\mj_A\:\kon\Sd_{p,B}$\> = \>
$\Sd_{p,A\kon B}\!\cirk(\mj_{p\kon\neg
p}\:\dis\c{\kon}_{B,A})\cirk d^R_{p\kon\neg
p,B,A}\cirk\!\c{\kon}_{A,(p\kon\neg p)\dis B}$,
\\*[1ex]
\` by $(\c{\kon}\c{\kon})$, $\mbox{($\c{\kon}$ {\it nat})}$,
$\mbox{($\Sd\kon\:\mj$)}$ and $\mbox{($\Sd$ {\it nat})}$.
\end{tabbing}

\noindent Note that in the stem-increasing equations the stem
index $B$ of $\Dk$ and $\Sd$ becomes more complex on the
right-hand sides, whereas the crown index $p$ does not change. We
have analogous stem-increasing equations for $\Sk$, $\Dkp\!$,
$\Skp\!$, $\Dd$, $\Sdp$ and $\Ddp\!$.

We will next prove several lemmata concerning \PN, which we will
find useful for calculations later on. For these lemmata we need
the following.

Let $\mds^{\neg p}$\index{DSnegp@$\mds^{\neg p}$ category} be the
category defined as \mds\ save that it is generated not by \Pe,
but by ${\Pe\cup\Pe^{\neg}}$. So the objects of $\mds^{\neg p}$
are formulae of ${\eL_{\kon,\dis}^{\neg p}}$, i.e.\ the objects of
\PN. For $A$ and $B$ formulae of ${\eL_{\kon,\dis}^{\neg p}}$, we
define when an occurrence of $p$ in $A$ is linked to an occurrence
of $p$ in $B$ in an arrow ${f\!:A\vdash B}$ of $\mds^{\neg p}$
analogously to what we had at the beginning of the preceding
section.

Let $\Xx$\index{Xixi@$\Xx$ arrow}\index{Xicon@$\Xk$
arrow}\index{Xidis@$\Xd$ arrow} for ${\!\ks\!\in\{\kon,\dis\}}$
stand for either $\Dx$, or $\Dxp\!$, or $\Sx$, or $\Sxp\!$, and
let a $\Xx_{B,A}$\emph{-term}\index{Xixiterm@$\Xx_{B,A}$-term} be
defined as a $\beta$-term in \S 3, save that $\beta$ is replaced
by $\Xx_{B,A}$. We use also $\Theta$\index{taulTheta@$\Theta$
arrow} as a variable alternative to $\Xi$. Then we have the
following.

\prop{$\Xk$-Permutation Lemma}{Let ${g\!:C\vdash D}$ be a
$\Xk_{p,B}$-term of \PN\ such that $x_1$ and $\neg x_2$ are
respectively the occurrences within $D$ of $p$ and $\neg p$ in the
crown of the head $\Xk_{p,B}$ of $g$, and let ${f\!:D\vdash E}$ be
an arrow term of $\mds^{\neg p}$ such that we have an occurrence
$y_1$ of $p$ and an occurrence $\neg y_2$ of $\neg p$ within a
subformula of $E$ of the form ${y_1\dis \neg y_2}$ or ${\neg
y_2\dis y_1}$, and $x_i$ is linked to $y_i$ for ${i\in\{1,2\}}$ in
$f$. Then there is a $\Tk_{p,B'}$-term ${g'\!:D'\vdash E}$ of \PN\
the crown of whose head is ${y_1\dis\neg y_2}$ or ${\neg y_2\dis
y_1}$, and there is an arrow term ${f'\!:C\vdash D'}$ of
$\mds^{\neg p}$ such that in \PN\ we have ${f\cirk g=g'\cirk
f'}$.}\index{Xicon-Permutation Lemma@$\Xk$-Permutation Lemma}

\dkz By the Development Lemma we can assume that $f$ is a
developed arrow term, and then it is enough to consider the case
when $f$ is either a $\beta$-term for $\beta$ a primitive arrow
term of $\mds^{\neg p}$ or $f$ is $\mj_E$. Note that in the
developed arrow term ${f_n\cirk\ldots\cirk f_1}$, which is equal
to $f$, we have that $f_1$ is $\mj_D$, and that $f_2$, if it
exists, cannot be a $d_{B,p,\neg p}$-term or a $d_{B,\neg
p,p}$-term such that $x_1$ and $\neg x_2$ are the occurrences of
$p$ and $\neg p$ in the right conjunct of the source ${B\kon(\neg
p\dis p)}$ or ${B\kon(p\dis\neg p)}$ of the head of $f_2$.
Otherwise, in the target of the head of $f_2$ we would obtain as
the left disjunct ${B\kon\neg p}$ or ${B\kon p}$, which together
with Lemma~2 would contradict the conditions put on $f$, and hence
also on ${f_n\cirk\ldots\cirk f_1}$, in the formulation of the
$\Xk$-Permutation Lemma.

The case when $f$ is $\mj_E$ is trivial, and there are also many
easy cases settled by bifunctorial and naturality equations. The
remaining, more interesting, cases are settled by the following
equations of \PN:

\begin{tabbing}
\mbox{\hspace{7em}}\= $\b{\kon}{\str}_{A,B,\neg p\dis
p}\!\cirk(\mj_A\:\kon \Dk_{p,B})=\;\Dk_{p,A\kon B}$, \` by
$\mbox{($\b{\kon}{}\Dk$)}$,\hspace{1em}
\\*[1ex]
\> $\b{\kon}{\rts}_{B_1,B_2,\neg p\dis p}\!\cirk\!\Dk_{p,B_1\kon
B_2}\;= \mj_{B_1}\:\kon\Dk_{p,B_2}$, \` by
$\mbox{($\b{\kon}{}\Dk$)}$,\hspace{1em}
\\[1.5ex]
\> $\b{\kon}{\str}_{A,\neg p\dis
p,B}\!\cirk(\mj_A\:\kon\Sk_{p,B})\;$\=$= \;\Dk_{p,A}\kon\:\mj_B$,
\` by $\mbox{($\b{\kon}{}\Dk\Sk$)}$,\hspace{1em}
\\*[1ex]
\> $\b{\kon}{\rts}_{B,\neg p\dis
p,A}\!\cirk(\Dk_{p,B}\kon\:\mj_A)$\>$= \mj_B\:\kon\Sk_{p,A}$, \`
by $\mbox{($\b{\kon}{}\Dk\Sk$)}$,\hspace{1em}
\\[1.5ex]
\> $\b{\kon}{\str}_{\neg p\dis p,B_1,B_2}\!\cirk\!\Sk_{p,B_1\kon
B_2}\;=\; \Sk_{p,B_1}\kon\:\mj_{B_2}$, \` by
$\mbox{($\b{\kon}{}\Sk$)}$,\hspace{1em}
\\*[1ex]
\> $\b{\kon}{\rts}_{\neg p\dis
p,B,A}\!\cirk(\Sk_{p,B}\kon\:\mj_A)=\; \Sk_{p,B\kon A}$, \` by
$\mbox{($\b{\kon}{}\Sk$)}$,\hspace{1em}
\\[1.5ex]
\> $\c{\kon}_{B,\neg p\dis p}\!\cirk\!\Dk_{p,B}\;=\; \Sk_{p,B}$,
\` by definition,\hspace{1em}
\\*[1ex]
\> $\c{\kon}_{\neg p\dis p,B}\!\cirk\!\Sk_{p,B}\;=\;\Dk_{p,B}$, \`
by definition and $(\c{\kon}\c{\kon})$,\hspace{1em}
\\[1.5ex]
\> $(\mj_B\:\kon\c{\dis}_{p,\neg
p})\cirk\!\Dk_{p,B}\;$\=$=\;\Dkp_{p,B}$, \` by
definition,\hspace{1em}
\\*[1ex]
\> $(\c{\dis}_{p,\neg
p}\kon\:\mj_B)\cirk\!\Sk_{p,B}$\>$=\;\Skp_{p,B}$, \` by definition
and $\mbox{($\c{\kon}$ {\it nat})}$,\hspace{1em}
\\[1.5ex]
\> $d_{\neg p\dis p,B_1,B_2}\cirk\!\Sk_{p,B_1\dis B_2}\;=\;
\Sk_{p,B_1}\dis\:\mj_{B_2}$, \` by
$\mbox{($d\!\Sk$)}$.\hspace{1em}
\end{tabbing}

\noindent Besides these equations, we have analogous equations
where ${\neg p\dis p}$ is replaced by ${p\dis \neg p}$, while
$\Dk$ and $\Sk$ are replaced by $\Dkp$ and $\Skp$ respectively,
and vice versa. \mbox{\hspace{1em}} \qed

\vspace{2ex}

We prove analogously the following dual of the preceding lemma.

\prop{$\Xd$-Permutation Lemma}{Let ${g\!:D\vdash C}$ be a
$\Xd_{p,B}$-term of \PN\ such that $x_1$ and $\neg x_2$ are
respectively the occurrences within $D$ of $p$ and $\neg p$ in the
crown of the head $\Xd_{p,B}$ of $g$, and let ${f\!:E\vdash D}$ be
an arrow term of $\mds^{\neg p}$ such that we have an occurrence
$y_1$ of $p$ and an occurrence $\neg y_2$ of $\neg p$ within a
subformula of $E$ of the form ${y_1\kon \neg y_2}$ or ${\neg
y_2\kon y_1}$, and $y_i$ is linked to $x_i$ for ${i\in\{1,2\}}$ in
$f$. Then there is a $\Td_{p,B'}$-term ${g'\!:E\vdash D'}$ of \PN\
the crown of whose head is ${y_1\kon\neg y_2}$ or ${\neg y_2\kon
y_1}$, and there is an arrow term ${f'\!:D'\vdash C}$ of
$\mds^{\neg p}$ such that in \PN\ we have ${g\cirk f=f'\cirk
g'}$.}\index{Xidis-Permutation Lemma@$\Xd$-Permutation Lemma}

Next we prove the following lemma, which involves the \emph{p-q-r}
Lemma of the preceding section.

\prop{\emph{p-$\neg$p-p} Lemma}{Let $x_1$, $\neg x_2$ and $x_3$ be
occurrences of $p$, $\neg p$ and $p$, respectively, in a formula
$A$ of ${\eL_{\kon,\dis}^{\neg p}}$, and let $y_1$, $\neg y_2$ and
$y_3$ be occurrences of $p$, $\neg p$ and $p$, respectively in a
formula $B$ of ${\eL_{\kon,\dis}^{\neg p}}$. Let ${\neg x_2\dis
x_3}$ or ${x_3\dis\neg x_2}$ be a subformula of $A$ and
${y_1\kon\neg y_2}$ or ${\neg y_2\kon y_1}$ a subformula of $B$.
Let ${g_1\!:A'\vdash A}$ be a $\Xk_{p,C}$-term of \PN\ such that
${\neg x_2\dis x_3}$ or ${x_3\dis\neg x_2}$ is the crown of the
head of $g_1$, let ${g_2\!:B\vdash B'}$ be a $\Td_{p,D}$-term of
\PN\ such that ${y_1\kon\neg y_2}$ or ${\neg y_2\kon y_1}$ is the
crown of the head of $g_2$, and let ${f\!:A\vdash B}$ be an arrow
term of $\mds^{\neg p}$ such that $x_i$ and $y_i$ are linked in
$f$ for ${i\in\{1,2,3\}}$. Then ${g_2\cirk f\cirk g_1}$ is equal
in \PN\ to an arrow term of $\mds^{\neg p}$.}\index{pc-negp-p
Lemma@\emph{p-$\neg$p-p} Lemma}

\dkz By the \emph{p-q-r} Lemma, ${f\!:A\vdash B}$ is equal in
$\mds^{\neg p}$, and hence also in \PN, to an arrow term of the
form ${f_y\cirk h\cirk f_x}$, where $h$ is a $d_{p,\neg
p,p}$-term, and the other conditions of the \emph{p-q-r} Lemma are
satisfied. So in \PN\ we have

\vspace{-.5ex}

\[g_2\cirk f\cirk g_1=g_2\cirk
f_y\cirk h\cirk f_x\cirk g_1=f_y'\cirk g_2'\cirk h\cirk g_1'\cirk
f_x',
\]

\vspace{-.5ex}

\noindent by the $\Xx$-Permutation Lemmata above. Here the head of
$g_1'$ must be $\Dk_{p,p}:p\vdash p\kon(\neg p\dis p)$, the head
of $h$ is ${d_{p,\neg p, p}\!:p\kon(\neg p\dis p)\vdash(p\kon\neg
p)\dis p}$, and the head of $g_2'$ must be ${\Sd_{p,p}:(p\kon\neg
p)\dis p\vdash p}$. By applying $\mbox{($\Sd\Dk$)}$, and perhaps
bifunctorial equations, we obtain that ${g_2'\cirk h\cirk g_1'}$
is equal in \PN\ to an arrow term of the form $\mj_A$, and hence
we have ${g_2\cirk f\cirk g_1=f_y'\cirk f_x'}$ in \PN, which
proves the lemma. \qed

\vspace{2ex}

To give an example of the application of the \emph{p-$\neg$p-p}
Lemma, consider the diagram in Figure 1.
\begin{table}
\begin{center}
\begin{picture}(200,340)

\put(117,5){\makebox(0,0)[r]{$p\kon q$}}
\put(120,35){\makebox(0,0)[r]{$(q\kon\neg q)\dis(p\kon q)$}}
\put(120,65){\makebox(0,0)[r]{$(q\kon(\!(p\,\kon\neg p)\dis\neg
q))\dis(p\kon q)$}}
\put(120,95){\makebox(0,0)[r]{$(q\kon(p\kon(\neg p\dis\:\neg
q)\!)\!)\dis(p\kon q)$}}
\put(120,125){\makebox(0,0)[r]{$((q\,\kon\, p) \kon(\neg p\dis\neg
q))\dis(p\kon q)$}} \put(120,155){\makebox(0,0)[r]{$((p\kon q)
\kon(\neg p\dis\neg q))\dis(p\kon q)$}}
\put(123,185){\makebox(0,0)[r]{$(p\kon q) \kon((\neg p\dis\!\neg
q)\,\dis(p\kon q))$}} \put(123,215){\makebox(0,0)[r]{$(p\kon q)
\kon((\neg q\dis\neg p)\dis(p\kon q))$}}
\put(128,245){\makebox(0,0)[r]{$(p\,\kon\, q) \kon(\neg q\dis(\neg
p\:\dis(p\kon q)))$}} \put(126,275){\makebox(0,0)[r]{$(p\kon
q)\:\kon(\neg q\!\dis((\neg p\:\dis\: p)\!\kon\! q))$}}
\put(70,305){\makebox(0,0)[r]{$(p\,\kon\, q) \kon(\neg q\dis q)$}}
\put(22,335){\makebox(0,0)[r]{$p\,\kon\, q$}}

\put(98,10){\line(0,1){19}} \put(114,10){\line(0,1){19}}
\put(98,40){\line(0,1){19}} \put(114,40){\line(0,1){19}}
\put(98,70){\line(0,1){19}} \put(114,70){\line(0,1){19}}
\put(98,100){\line(0,1){19}} \put(114,100){\line(0,1){19}}
\put(98,130){\line(0,1){19}} \put(114,130){\line(0,1){19}}
\put(98,160){\line(0,1){19}} \put(114,160){\line(0,1){19}}
\put(98,190){\line(0,1){19}} \put(114,190){\line(0,1){19}}
\put(98,220){\line(0,1){19}} \put(114,220){\line(0,1){19}}
\put(98,250){\line(0,1){19}} \put(114,250){\line(0,1){19}}

\put(71,40){\line(0,1){19}} \put(71,70){\line(0,1){19}}
\put(71,100){\line(0,1){19}} \put(71,130){\line(0,1){19}}
\put(71,160){\line(0,1){19}}

\put(46,70){\line(0,1){19}} \put(46,100){\line(0,1){19}}
\put(46,130){\line(0,1){19}} \put(46,160){\line(0,1){19}}

\put(18,70){\line(0,1){19}} \put(18,100){\line(0,1){19}}

\put(-3,70){\line(0,1){19}} \put(-3,100){\line(0,1){19}}

\put(0,160){\line(0,1){19}} \put(0,190){\line(0,1){19}}
\put(0,220){\line(0,1){19}} \put(0,250){\line(0,1){19}}
\put(0,280){\line(0,1){19}} \put(0,310){\line(0,1){20}}

\put(18,160){\line(0,1){19}} \put(18,190){\line(0,1){19}}
\put(18,220){\line(0,1){19}} \put(18,250){\line(0,1){19}}
\put(18,280){\line(0,1){19}} \put(18,310){\line(0,1){20}}

\put(-3,130){\line(1,1){19}} \put(18,130){\line(-1,1){19}}
\put(46,190){\line(1,1){20}} \put(71,190){\line(-1,1){19}}

\put(48,220){\line(0,1){19}} \put(73,220){\line(0,1){19}}
\put(48,250){\line(0,1){19}} \put(73,250){\line(0,1){19}}
\put(48,280){\line(0,1){19}} \put(114,280){\line(-5,2){49}}

\put(47,40){\line(-5,2){49}}

\put(62,29){\oval(20,20)[b]} \put(31,59){\oval(24,10)[b]}
\put(55,310){\oval(16,16)[t]} \put(86,280){\oval(24,10)[t]}

\put(130,20){\makebox(0,0)[l]{$\Sd_{q,p\kon q}$}}

\put(130,50){\makebox(0,0)[l]{$(\mj_q\:\kon\Sd_{p,\neg
q})\dis\mj_{p\kon q}$}}

\put(130,80){\makebox(0,0)[l]{$(\mj_q\kon d_{p,\neg p,\neg
q})\dis\mj_{p\kon q}$}}

\put(130,110){\makebox(0,0)[l]{$\b{\kon}{\rts}_{q,p,\neg p\dis\neg
q}\dis\:\mj_{p\kon q}$}}

\put(130,140){\makebox(0,0)[l]{$(\c{\kon}_{p,q}\kon\:\mj_{\neg
p\dis\neg q})\dis\mj_{p\kon q}$}}

\put(130,170){\makebox(0,0)[l]{$d_{p\kon q,\neg p\dis\neg q,p\kon
q}$}}

\put(130,200){\makebox(0,0)[l]{$\mj_{p\kon q}\kon(\c{\dis}_{\neg
p,\neg q}\dis\:\mj_{p\kon q})$}}

\put(130,230){\makebox(0,0)[l]{$\mj_{p\kon
q}\:\kon\b{\dis}{\str}_{\neg q,\neg p,p\kon q}$}}

\put(130,260){\makebox(0,0)[l]{$\mj_{p\kon q}\kon(\mj_{\neg q}\dis
d^R_{\neg p,p,q})$}}

\put(130,295){\makebox(0,0)[l]{$\mj_{p\kon q}\kon(\mj_{\neg
q}\:\dis\Sk_{p,q})$}}

\put(130,325){\makebox(0,0)[l]{$\Dk_{q,p\kon q}$}}

\end{picture}

\vspace{6ex}

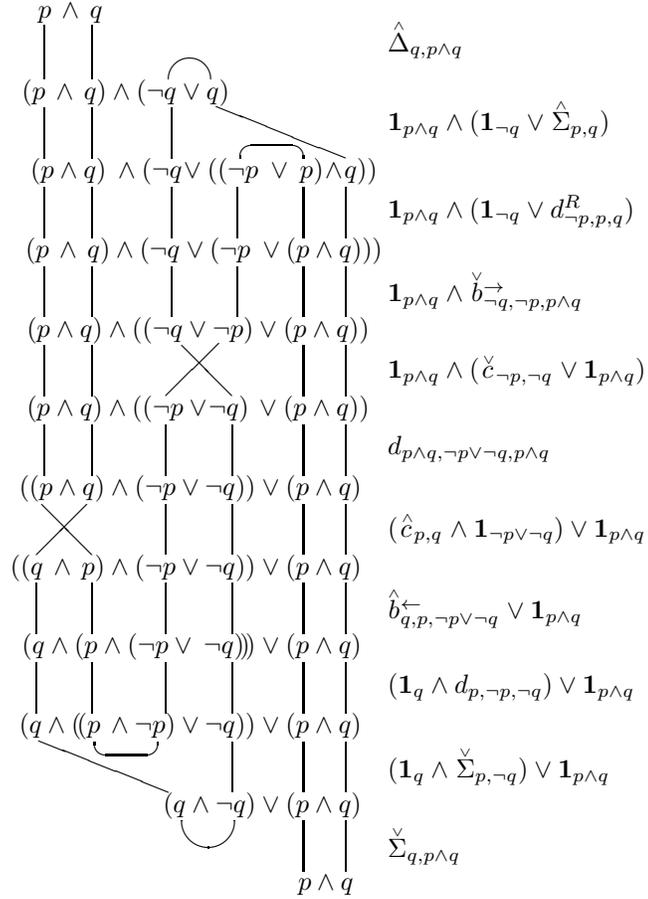
{\sc Figure~1}
\end{center}
\end{table}
\noindent This diagram corresponds to ${G(\Sd_{q,p\kon q}\!\cirk
h\cirk\!\Dk_{q,p\kon q})}$ for an arrow term $h$ of \PN, which is
of the form ${g_2\cirk f\cirk g_1}$ for $g_1$ being ${\mj_{p\kon
q}\kon(\mj_{\neg q}\:\dis\Sk_{p,q})}$, $g_2$ being
${(\mj_q\kon\Sd_{p,\neg q})\dis\mj_{p\kon q}}$ and $f$ an arrow
term of $\mds^{\neg p}$. Then by applying the
\mbox{\emph{p-$\neg$p-p}} Lemma we obtain an arrow term $f'$ of
$\mds^{\neg p}$ equal to ${g_2\cirk f\cirk g_1}$ in \PN, and next
by applying the \emph{p-$\neg$p-p} Lemma (as a matter of fact, the
$q$-$\neg q$-$q$ Lemma), we obtain an arrow term $h'$ of
$\mds^{\neg p}$ equal to ${\Sd_{q,p\kon q}\!\cirk
f'\cirk\!\Dk_{q,p\kon q}}$ in \PN. By \mds\ Coherence of \S 4, we
may conclude that $h'$, and hence also ${\Sd_{q,p\kon q}\!\cirk
h\cirk\!\Dk_{q,p\kon q}}$, is equal to $\mj_{p\kon q}$ in \PN.

Here is a lemma analogous to the \emph{p-$\neg$p-p} Lemma.

\prop{\emph{$\neg$p-p-$\neg$p} Lemma}{Let $\neg x_1$, $x_2$ and
$\neg x_3$ be occurrences of $\neg p$, $p$ and $\neg p$,
respectively, in a formula $A$ of ${\eL_{\kon,\dis}^{\neg p}}$,
and let $\neg y_1$, $y_2$ and $\neg y_3$ be occurrences of $\neg
p$, $p$ and $\neg p$, respectively, in a formula $B$ of
${\eL_{\kon,\dis}^{\neg p}}$. Let ${g_1\!:A'\vdash A}$ be a
$\Xk_{p,C}$-term of \PN\ such that ${x_2\dis\neg x_3}$ or ${\neg
x_3\dis x_2}$ is the crown of the head of $g_1$, let
${g_2\!:B\vdash B'}$ be a $\Td_{p,D}$-term of \PN\ such that
${\neg y_1\kon y_2}$ or ${y_2\kon \neg y_1}$ is the crown of the
head of $g_2$, and let ${f\!:A\vdash B}$ be an arrow term of
$\mds^{\neg p}$ such that $x_i$ and $y_i$ are linked in $f$ for
${i\in\{1,2,3\}}$. Then ${g_2\cirk f\cirk g_1}$ is equal in \PN\
to an arrow term of $\mds^{\neg p}$.}\index{pc-negp-pp
Lemma@\emph{$\neg$p-p-$\neg$p} Lemma}

\noindent To prove this lemma we proceed as for the
\emph{p-$\neg$p-p} Lemma, relying on the equation
$\mbox{($\Sdp\Dkp$)}$ of \PN.

\section{\large\bf The equivalence of \PNN\ and \PN}

In this section we show that the categories \PNN\ and \PN\ are
equivalent categories. We define inductively a functor
$F$\index{F@$F$ functor} from the category \PNN\ to \PN\ in the
following manner. On objects we have

\begin{tabbing}
\mbox{\hspace{2em}}\= $Fp=p$,\quad for $p$ a letter,
\\[1ex]
\> $F(A\ks B)=FA\ks FB$,\quad for $\!\ks\!\in\{\kon,\dis\}$,
\\[1ex]
\> $F\neg p=\neg p$,\quad for $p$ a letter,
\\[1ex]
\> $F\neg\neg A=FA$,
\\[1ex]
\> $F\neg(A\kon B)=F\neg A\dis F\neg B$,
\\*[1ex]
\> $F\neg(A\dis B)=F\neg A\kon F\neg B$.
\end{tabbing}

\noindent On arrows we have

\begin{tabbing}
\mbox{\hspace{2em}}\=
$F\alpha_{A_1,\ldots,A_n}=\alpha_{FA_1,\ldots,FA_n}$,
\\*[.5ex]
for $\alpha_{A_1,\ldots,A_n}$ being $\mj_A$,
$\b{\xi}{\str}_{A,B,C}$, $\b{\xi}{\rts}_{A,B,C}$, $\c{\xi}_{A,B}$
or $d_{A,B,C}$ where $\!\ks\!\in\{\kon,\dis\}$,
\\[2ex]
\> $F\!\Dk_{p,A}\;$\=$=\;\Dk_{p,FA}:FA\vdash FA\kon(\neg p\dis
p)$,
\\*[1ex]
\> $F\!\Sd_{p,A}$\>$=\;\Sd_{p,FA}:(p\kon\neg p)\dis FA\vdash FA$,
\\[2ex]
\> $F\!\Dk_{\neg B,A}$ \= = \= $(\mj_{FA}\:\kon\c{\dis}_{FB,F\neg
B})\cirk F\!\Dk_{B,A}:FA\vdash FA\kon(FB\dis F\neg B)$,
\\*[1ex]
\> $F\!\Sd_{\neg B,A}$ \> = \>
$F\!\Sd_{B,A}\!\cirk(\c{\kon}_{F\neg B,FB}\dis\:\mj_{FA})\!:(F\neg
B\kon FB)\dis FA\vdash FA$,
\\[2ex]
\pushtabs \mbox{\hspace{2em}}\= $F\!\Dk_{B\kon C,A}$ \= = \=
$(\mj_{FA}\:\kon((\c{\dis}_{F\neg B,F\neg C}\dis\:\mj_{FB\kon
FC})\cirk\!\b{\dis}{\str}_{F\neg C,F\neg B,FB\kon FC}\!\cirk$
\\*[.5ex]
\` $\cirk(\mj_{F\neg C}\dis(d^R_{F\neg B,FB,FC}\cirk
\c{\kon}_{FC,F\neg B\dis FB}\!\cirk F\!\Dk_{B,C}))))\cirk
F\!\Dk_{C,A}:$
\\*[.5ex]
\` $FA\vdash FA\kon((F\neg B\dis F\neg C)\dis (FB\kon FC))$,
\\*[1ex]
 \> $F\!\Sd_{B\kon C,A}$ \> = \>
 $F\!\Sd_{C,A}\!\cirk((\mj_{FC}\kon(F\!\Sd_{B,\neg C}\!\cirk
d_{FB,F\neg B,F\neg C}))\cirk$
\\*[.5ex]
\` $\cirk\!\b{\kon}{\rts}_{FC,FB,F\neg B\dis F\neg
C}\!\cirk(\c{\kon}_{FB,FC}\kon\mj_{F\neg B\dis F\neg
C}))\dis\mj_{FA})\!:$
\\*[.5ex]
\` $((FB\kon FC)\kon(F\neg B\dis F\neg C))\dis FA\vdash FA$,
\\[2ex]
\> $F\!\Dk_{B\dis C,A}$ \> = \> $(\mj_{FA}\kon((\c{\kon}_{F\neg
C,F\neg B}\dis\:\mj_{FB\dis FC})\cirk\!\b{\dis}{\rts}_{F\neg C\kon
F\neg B,FB,FC}\!\cirk$
\\*[.5ex]
\` $\cirk((d_{F\neg C,F\neg B,FB}\cirk F\!\Dk_{B,\neg
C})\dis\mj_{FC})))\cirk F\!\Dk_{C,A}:$
\\*[.5ex]
\` $FA\vdash FA\kon((F\neg B\kon F\neg C)\dis (FB\dis FC))$,
\\*[1ex]
\> $F\!\Sd_{B\dis C,A}$ \> = \>
$F\!\Sd_{C,A}\!\!\cirk(((F\!\Sd_{B,C}\!\!\cirk\!\!\c{\dis}_{FB\kon
F\neg B,FC}\!\!\cirk d^R_{FC,FB,F\neg B})\!\kon\!\mj_{F\neg
C})\!\cirk$
\\*[.5ex]
\` $\cirk\!\b{\kon}{\str}_{FC\dis FB,F\neg B,F\neg
C}\!\!\cirk(\c{\dis}_{FC,FB}\kon\:\mj_{F\neg B\kon F\neg
C}))\dis\mj_{FA})\!:$
\\*[.5ex]
\` $((FB\dis FC)\kon(F\neg B\kon F\neg C))\dis FA\vdash FA$,
\\[1.5ex]
\pushtabs \mbox{\hspace{2em}}\= $F(f\cirk g)$ \= = \= $Ff\cirk
Fg$,
\\[1ex]
\> $F(f\ks g)$ \> = \> $Ff\ks Fg$, \quad for
$\!\ks\!\in\{\kon,\dis\}$. \poptabs \poptabs
\end{tabbing}

It is easy to infer

\begin{tabbing}
\mbox{\hspace{2em}}\= $F\!\Dk_{\neg B,A}$ \= = \=
$F\!\Dkp_{B,A}$,\hspace{6em}\= $F\!\Sd_{\neg B,A}$ \= = \=
$F\!\Sdp_{B,A}$,
\\[1ex]
\>$F\!\Dkp_{\neg B,A}$\>=\>$F\!\Dk_{B,A}$,\>$F\!\Sdp_{\neg
B,A}$\>=\>$F\!\Sd_{B,A}$,
\\[1ex]
\> $F\!\Dk_{B,A}$ \> = \> $F\!\Dk_{B,FA}$, \> $F\!\Sd_{B,A}$ \> =
\> $F\!\Sd_{B,FA}$.
\end{tabbing}

To ascertain that $F$ so defined is indeed a functor, we have to
verify that if ${f=g}$ is an instance of one of the \PN\
equations, then ${Ff=Fg}$ holds in \PN. This is done by induction
on the number od occurrences of connectives in the crown indices
occurring in these equations.

For $\mbox{($\Dk$ {\it nat})}$ and $\mbox{($\Sd$ {\it nat})}$ this
is a very easy matter. For $\mbox{($\b{\kon}{}\Dk$)}$,
$\mbox{($\b{\dis}{}\Sd$)}$, $\mbox{($d\!\Sk$)}$ and
$\mbox{($d\!\Dd$)}$ we use essentially naturality equations. (In
that context, it might be easier to rely on the equations
$\mbox{($d^R\!\Dk$)}$ and $\mbox{($d^R\!\Sd$)}$, which are
alternative to $\mbox{($d\!\Sk$)}$ and $\mbox{($d\!\Dd$)}$.)

To verify $\mbox{($\Sd\Dk$)}$ in cases where $A$ is of the form
${B\kon C}$ or ${B\dis C}$, we rely on the induction hypothesis
that if ${f=g}$ is an instance of a \PN\ equation such that the
crown indices are $B$ and $C$, then we have ${Ff=Fg}$ in \PN. This
induction hypothesis entails that we can proceed as in the proof
of the \emph{p-$\neg$p-p} Lemma in the preceding section, first
for $p$ replaced by $B$, and then for $p$ replaced by $C$.
Finally, we apply \mds\ Coherence (see the example at the end of
the preceding section). To verify $\mbox{($\Sd\Dk$)}$ in case $A$
is of the form $\neg B$, we rely on the induction hypothesis for
the equation $\mbox{($\Sdp\Dkp$)}$.

To verify $\mbox{($\Sdp\Dkp$)}$ we proceed analogously. In case
$A$ is ${B\kon C}$ or ${B\dis C}$, we rely on the proof of the
\emph{$\neg$p-p-$\neg$p} Lemma in the preceding section, and in
case $A$ is $\neg B$ we rely on the induction hypothesis for the
equation $\mbox{($\Sd\Dk$)}$. This concludes the verification that
$F$ is a functor from \PNN\ to \PN.

In the definition of $F$, there is some freedom in choosing the
clauses for ${F\!\Xx_{B\psi C,A}}$, where
${\Xi\in\{\Delta,\Sigma\}}$ and ${\!\ks\!,
\mbox{\footnotesize$\psi$}\in\{\kon,\dis\}}$. Ours enable us to
apply easily the \emph{p-$\neg$p-p} and \emph{$\neg$p-p-$\neg$p}
Lemmata in verifying that $F$ is a functor.

We define a functor $\Fn$\index{Faaneg@$\Fn$ functor} from \PN\ to
\PNN\ by stipulating that $\Fn A=A$ and $\Fn f=f$. It is clear
that if ${f=g}$ in \PN, then ${\Fn f=\Fn g}$ in \PNN; so $\Fn$ is
indeed a functor.

Our purpose is to show that \PNN\ and \PN\ are equivalent
categories via the functors $F$ and $\Fn$. It is clear that ${F\Fn
A=A}$ and ${F\Fn f=f}$. Since ${\Fn FA=FA}$, we have to define in
\PNN\ an isomorphism ${i_A\!:A\vdash FA}$. For that we need the
following auxiliary definitions in \PNN:

\begin{tabbing}
\mbox{\hspace{0pt}}\= $\r{\kon}{\str}_{A,B}$ \= $=_{\df}$ \=\kill
\> \> $\n{\str}_A\!$\' $=_{\df}$ \> $\Sdp_{\neg A,A}\!\cirk
d_{\neg\neg A,\neg A,A}\cirk\!\Dk_{A,\neg\neg A}$ \= : \=
$\neg\neg A\vdash A$,\index{na@$\n{\str}$ arrow}
\\*[1ex]
\> \> $\n{\rts}_A\!$\' $=_{\df}$ \> $\Sd_{A,\neg\neg A}\!\cirk
d_{A,\neg A,\neg\neg A}\cirk\!\Dkp_{\neg A,A}$\> :\>
$A\vdash\neg\neg A$,\index{nb@$\n{\rts}$ arrow}
\\[1.5ex]
\> $\r{\kon}{\str}_{A,B}$ \> $=_{\df}$ \> $\Sdp_{A\kon B,\neg
A\dis\neg B}\!\!\cirk d_{\neg(A\kon B),A\kon B,\neg A\dis\neg
B}\!\cirk(\mj_{\neg(A\kon B)}\kon((\mj_{A\kon
B}\dis\!\c{\dis}_{\neg A,\neg B})\!\cirk$
\\*[.5ex]
\` $\cirk\!\b{\dis}{\rts}_{A\kon B,\neg B,\neg
A}\!\!\cirk((d_{A,B,\neg B}\cirk\!\Dkp_{B,A})\dis\mj_{\neg
A})))\cirk\!\Dkp_{A,\neg(A\kon B)}\::$
\\*[.5ex]
\` $\neg(A\kon B)\vdash\neg A\dis\neg
B$,\index{rcona@$\r{\kon}{\str}$ arrow}
\\[1ex]
\> $\r{\kon}{\rts}_{A,B}$ \> $=_{\df}$ \> $\Sdp_{A,\neg(A\kon
B)}\!\cirk((((\Ddp_{B,\neg A}\!\cirk d^R_{\neg A,\neg
B,B})\kon\mj_A)\cirk\!\b{\kon}{\str}_{\neg A\dis\neg
B,B,A}\!\!\cirk$
\\*[.5ex]
\` $\cirk(\mj_{\neg A\dis\neg
B}\,\kon\c{\kon}_{A,B}))\dis\mj_{\neg(A\kon B)})\cirk d_{\neg
A\dis\neg B,A\kon B,\neg(A\kon B)}\cirk\!\Dkp_{A\kon B,\neg
A\dis\neg B}\::$
\\*[.5ex]
\` $\neg A\dis\neg B\vdash\neg(A\kon
B)$,\index{rconb@$\r{\kon}{\rts}$ arrow}
\\[1.5ex]
\> $\r{\dis}{\str}_{A,B}$ \> $=_{\df}$ \> $\Sdp_{A\dis B,\neg
A\kon\neg B}\!\!\cirk d_{\neg(A\dis B),A\dis B,\neg A\kon\neg
B}\cirk(\mj_{\neg(A\dis B)}\!\kon\!((\c{\dis}_{A,B}\!\dis\mj_{\neg
A\kon\neg B})\!\cirk$
\\*[.5ex]
\` $\cirk\!\b{\dis}{\str}_{B,A,\neg A\kon \neg
B}\!\!\cirk(\mj_B\dis( d^R_{A,\neg A,\neg B}\cirk\!\Skp_{A,\neg
B}))))\cirk\!\Dkp_{B,\neg(A\dis B)}\::$
\\*[.5ex]
\` $\neg(A\dis B)\vdash \neg A\kon\neg
B$,\index{rdisa@$\r{\dis}{\str}$ arrow}
\\*[1ex]
\> $\r{\dis}{\rts}_{A,B}$ \> $=_{\df}$ \> $\Sdp_{B,\neg(A\dis
B)}\!\cirk(((\mj_{\neg B}\kon(\Sdp_{A,B}\!\cirk d_{\neg
A,A,B}))\cirk\!\b{\kon}{\rts}_{\neg B,\neg A,A\dis B}\!\!\cirk$
\\*[.5ex]
\` $\cirk(\c{\kon}_{\neg A,\neg B}\kon\:\mj_{A\dis
B}))\dis\mj_{\neg(A\dis B)})\cirk d_{\neg A\kon \neg B,A\dis
B,\neg(A\dis B)}\cirk\!\Dkp_{A\dis B,\neg A\kon\neg B}\::$
\\*[.5ex]
\` $\neg A\kon\neg B\vdash\neg(A\dis
B)$.\index{rdisb@$\r{\dis}{\rts}$ arrow}
\end{tabbing}

It can be shown that in \PNN\ we have the following equations:

\begin{tabbing}
\mbox{\hspace{2em}}\=
$\r{\kon}{\str}_{A,B}\!\cirk\!\r{\kon}{\rts}_{A,B}\;$\=$=\mj_{\neg
A\dis\neg B}$,\hspace{6em}\=
$\r{\kon}{\rts}_{A,B}\!\cirk\!\r{\kon}{\str}_{A,B}\;$\=$=\mj_{\neg(A\kon
B)}$,\kill

\>\mbox{\hspace{1.5em}}$\n{\str}_A\cirk\n{\rts}_A$\>$=\mj_A$,\>
\mbox{\hspace{1.5em}}$\n{\rts}_A\cirk\n{\str}_A$\>$=\mj_{\neg\neg
A}$,
\\[1ex]
\>
$\r{\kon}{\str}_{A,B}\!\cirk\!\r{\kon}{\rts}_{A,B}$\>$=\mj_{\neg
A\dis\neg B}$,\>
$\r{\kon}{\rts}_{A,B}\!\cirk\!\r{\kon}{\str}_{A,B}$\>$=\mj_{\neg(A\kon
B)}$,
\\[1ex]
\>
$\r{\dis}{\str}_{A,B}\!\cirk\!\r{\dis}{\rts}_{A,B}$\>$=\mj_{\neg
A\kon\neg B}$,\>
$\r{\dis}{\rts}_{A,B}\!\cirk\!\r{\dis}{\str}_{A,B}$\>$=\mj_{\neg(A\dis
B)}$,
\end{tabbing}

\noindent which means that $\n{\str}$ and $\n{\rts}$, as well as
$\r{\xi}{\str}$ and $\r{\xi}{\rts}$ are inverses of each other. To
derive these equations in \PNN, we use essentially $\mbox{($\Dk$
{\it nat})}$, $\mbox{($\Sd$ {\it nat})}$, the \emph{p-$\neg$p-p}
and \emph{$\neg$p-p-$\neg$p} Lemmata, and \mds\ Coherence. (If an
equation holds in \PN, then every substitution instance of it
obtained by replacing letters uniformly by formulae of
$\eL_{\neg,\kon,\dis}$ holds in \PNN; this enables us to apply the
\emph{p-$\neg$p-p} and \emph{$\neg$p-p-$\neg$p} Lemmata.) The
definitions of $\n{\str}$, $\n{\rts}$, $\r{\xi}{\str}$ and
$\r{\xi}{\rts}$, for ${\!\ks\!\in\{\kon,\dis\}}$, are such that
they enable an easy application of the \emph{p-$\neg$p-p} and
\emph{$\neg$p-p-$\neg$p} Lemmata.

Then we define ${i_A\!:A\vdash FA}$\index{ia@$i_A$ arrow} and its
inverse ${i^{-1}_A\!:FA\vdash A}$\index{iaa@$i^{-1}_A$ arrow} by
induction on the complexity of the formula $A$ of
$\eL_{\neg,\kon,\dis}$ (cf.\ \cite{DP04}, Section 14.1):

\begin{tabbing}
\mbox{\hspace{1em}}\= $i_{\neg(A_1\kon A_2)}=(i_{\neg A_1}\dis
i_{\neg A_2})\cirk\!
\r{\kon}{\str}_{A_1,A_2}$,\mbox{\hspace{1em}}\=$i^{-1}_{\neg(A_1\kon
A_2)}=\:\r{\kon}{\rts}_{A_1,A_2}\!\!\cirk (i^{-1}_{\neg A_1}\dis
i^{-1}_{\neg A_2})$,\kill

\mbox{\hspace{10em}} $i_A=i^{-1}_A=\mj_A$, \quad if $A$ is $p$ or
$\neg p$, for $p$ a letter,
\\[1ex]
\> $i_{A_1\kst A_2}=i_{A_1}\ks i_{A_2}$, \> $i^{-1}_{A_1\kst
A_2}=i^{-1}_{A_1}\ks i^{-1}_{A_2}$, \` for
$\!\ks\!\in\{\kon,\dis\}$,
\\[1ex]
\> $i_{\neg\neg B}=i_B\cirk\n{\str}_B$,\> $i^{-1}_{\neg\neg
B}=\n{\rts}_B\cirk i^{-1}_B$,
\\[1ex]
\> $i_{\neg(A_1\kon A_2)}=(i_{\neg A_1}\dis i_{\neg A_2})\cirk\!
\r{\kon}{\str}_{A_1,A_2}$,\>$i^{-1}_{\neg(A_1\kon
A_2)}=\:\r{\kon}{\rts}_{A_1,A_2}\!\!\cirk (i^{-1}_{\neg A_1}\dis
i^{-1}_{\neg A_2})$,
\\[1ex]
\>$i_{\neg(A_1\dis A_2)}=(i_{\neg A_1}\kon i_{\neg A_2})\cirk\!
\r{\dis}{\str}_{A_1,A_2}$,\>$i^{-1}_{\neg(A_1\dis
A_2)}=\:\r{\dis}{\rts}_{A_1,A_2}\!\!\cirk (i^{-1}_{\neg A_1}\kon
i^{-1}_{\neg A_2})$.

\end{tabbing}

\noindent We can then prove the following (cf.\ \cite{DP04},
Section 14.1).

\prop{Auxiliary Lemma}{For every arrow term ${f\!:A\vdash B}$ of
\PNN\ we have $f=i^{-1}_B\cirk Ff\cirk i_A$ in \PNN.}

\dkz We proceed by induction on the complexity of the arrow term
$f$. If $f$ is a primitive arrow term $\mj_A$,
$\b{\xi}{\str}_{A,B,C}$, $\b{\xi}{\rts}_{A,B,C}$, $\c{\xi}_{A,B}$
or $d_{A,B,C}$, for ${\!\ks\!\in\{\kon,\dis\}}$, then we use
naturality equations, and the fact that $i_D$ is an isomorphism.

If $f$ is $\Dk_{D,A}$, then we proceed by induction on the
complexity of $D$. (This is an auxiliary induction in the basis of
the main induction.) If $D$ is $p$, then we use $\mbox{($\Dk$ {\it
nat})}$ and the fact that $i_A$ is an isomorphism.

If $D$ is $\neg B$, then we rely on the following equation of
\PNN:

\begin{tabbing}
\mbox{\hspace{2em}}\= $\mbox{($\Dk\! n$)}$\quad\quad $\Dk_{\neg
B,A}\;=(\mj_A\kon(\n{\rts}_B\dis\mj_{\neg
B}))\cirk\!\Dkp_{B,A}$,\index{Deltaconn@$\mbox{($\Dk n$)}$
equation}
\end{tabbing}

\noindent together with the induction hypothesis. To derive
$\mbox{($\Dk\! n$)}$ we have

\begin{tabbing}

$(\mj_A\kon(\n{\rts}_B\dis\mj_{\neg B}))\cirk\!\Dkp_{B,A}$
\\[1ex]
\mbox{\hspace{2em}} \= $=(\mj_A\kon(\Sd_{B,\neg\neg
B}\dis\:\mj_{\neg B}))\cirk(\mj_A\kon(d_{B,\neg B,\neg\neg
B}\dis\mj_{\neg B}))\cirk$
\\*[.5ex]
\` $\cirk(\mj_A\kon(\Dkp_{\neg B,B}\dis\:\mj_{\neg
B}))\cirk\!\Dkp_{B,A}$, \quad by bifunctorial equations,
\\[1ex]
\> $=(\mj_A\kon(\Sd_{B,\neg\neg B}\dis\mj_{\neg
B}))\cirk(\mj_A\kon((d_{B,\neg B,\neg\neg B}\dis\mj_{\neg
B})\cirk$
\\*[.5ex]
\` $\cirk\!\c{\dis}_{B\kon(\neg B\dis\neg\neg B),\neg B}\cirk
d^R_{\neg B, B,\neg B\dis\neg\neg B}\cirk(\c{\dis}_{\neg
B,B}\kon\:\mj_{\neg B\dis\neg\neg B})))\cirk$
\\*[.5ex]
\` $\cirk\!\b{\kon}{\rts}_{A,\neg B,B\dis\neg
B}\!\cirk(\Dkp_{B,A}\kon\:\mj_{\neg B\dis\neg\neg
B})\cirk(\mj_A\:\kon\c{\dis}_{\neg B,\neg\neg B})\cirk\!\Dk_{\neg
B,A}$,
\end{tabbing}

\vspace{-1ex}

\noindent by stem-increasing equations involving $\Dkp$ analogous
to $\mbox{($\mj\:\dis\Dk$)}$ and $\mbox{($\mj\:\kon\Dk$)}$ of the
preceding section, and also $\mbox{($\Dkp$ {\it nat})}$. The
equation $\mbox{($\Dk\! n$)}$ follows by applying the
\emph{$\neg$p-p-$\neg$p} Lemma (with $p$ replaced by $A$), and
\mds\ Coherence.

If $D$ is ${B\kon C}$, then we rely on the following equation of
\PNN:

\begin{tabbing}
\mbox{\hspace{2em}}\= $\mbox{($\Dk\! r$)}$\quad\quad $\Dk_{B\kon
C,A}\;=(\mj_A\kon(((\r{\kon}{\rts}_{B,C}\!\cirk\!\c{\dis}_{\neg
B,\neg C})\dis\mj_{B\kon C})\cirk\!\b{\dis}{\str}_{\neg C,\neg
B,B\kon C}\!\cirk$
\\
\` $\cirk(\mj_{\neg C}\dis(d^R_{\neg
B,B,C}\cirk\!\Sk_{B,C}))))\cirk\!\Dk_{C,A}$,\index{Deltaconnat1@$\mbox{($\Dk
r$)}$ equation}
\end{tabbing}

\noindent together with the induction hypothesis. To show that
$\mbox{($\Dk\! r$)}$ holds in \PNN\ we proceed as above, by
applying essentially stem-increasing equations together with the
\emph{p-$\neg$p-p} Lemma. We proceed analogously when $D$ is
${B\dis C}$.

The cases we have if $f$ is $\Sd_{D,A}$ are dual to those we had
above for $f$ being $\Dk_{D,A}$. In all these cases we proceed in
an analogous manner. This concludes the basis of the induction.

If $f$ is ${f_2\cirk f_1}$, then by the induction hypothesis we
have

\[
f_2\cirk f_1=i^{-1}_B\cirk Ff_2\cirk i_C\cirk i^{-1}_C\cirk
Ff_1\cirk i_A
\]
which yields ${f=i^{-1}_B\cirk Ff\cirk i_A}$, by the fact that
$i_C$ is an isomorphism and by the functoriality of $F$.

If $f$ is ${f_1\ks f_2}$, for ${\!\ks\!\in\{\kon,\dis\}}$, then
$i_{A_1\kst A_2}$ is ${i_{A_1}\ks i_{A_2}}$ and $i^{-1}_{B_1\kst
B_2}$ is ${i^{-1}_{B_1}\ks i^{-1}_{B_2}}$; we obtain
${f=i^{-1}_B\cirk Ff\cirk i_A}$ by using bifunctorial equations.
\qed

\vspace{2ex}

The Auxiliary Lemma shows that $i_A$ is an isomorphism natural in
$A$, and so we may conclude that \PNN\ and \PN\ are equivalent
categories.

\section{\large\bf \PN\ Coherence}

We define a functor $G$\index{G functor@$G$ functor} from \PN\ to
\emph{Br} as we defined it from \PNN\ to \emph{Br}. In the clauses
for $\Dk_{B,A}$ and $\Sd_{B,A}$ we just restrict $B$ to a letter
$p$. For $f$ an arrow term of \PNN\ we have that $GFf$ coincides
with $Gf$ where $F$ is the functor from \PNN\ to \PN\ of the
preceding section, $G$ in $GFf$ is the functor $G$ from \PN\ to
\emph{Br} and $G$ in $Gf$ is the functor $G$ from \PNN\ to
\emph{Br}. To show that, it is essential to check that
${GF\!\Dk_{B,A}}$ and ${GF\!\Sd_{B,A}}$ coincide with
${G\!\Dk_{B,A}}$ and ${G\!\Sd_{B,A}}$ respectively.

In this section we will prove that $G$ from \PN\ to \emph{Br} is
faithful. This will imply that $G$ from \PNN\ to \emph{Br} is
faithful too.

Analogously to what we had at the beginning of \S 5, we define
when an occurrence $x$ of a letter $p$ in $A$ is {\it linked} to
an occurrence $y$ of the same letter $p$ in $B$ in an arrow
${f\!:A\vdash B}$ of \PN. We say that $x$ and $y$ are {\it
directly linked} in a headed factorized arrow term
${f_n\cirk\ldots\cirk f_1}$ of \PN\ when $x$ and $y$ are linked in
the arrow ${f_n\cirk\ldots\cirk f_1}$, and for every
${i\in\{2,\ldots,n\}}$ if $f_i$ is a $\Sd_{p,C}$-term and $z$ is
one of the two occurrences of $p$ in the crown ${p\kon\neg p}$ of
the head of $f_i$, then $x$ and $z$ are not linked in the arrow
${f_{i-1}\cirk\ldots\cirk f_1}$ (see the end of \S 3 for the
definition of headed factorized arrow term).

An alternative definition of directly linked $x$ and $y$ in a
headed factorized arrow term ${f_1\cirk\ldots\cirk f_n}$ of \PN\
is obtained by stipulating that $x$ and $y$ are linked in the
arrow ${f_1\cirk\ldots\cirk f_n}$, and for every
${i\in\{2,\ldots,n\}}$ if $f_i$ is a $\Dk_{p,D}$-term and $z$ is
one of the two occurrences of $p$ in the crown ${\neg p\dis p}$ of
the head of $f_i$, then $z$ and $y$ are not linked in the arrow
${f_1\cirk\ldots\cirk f_{i-1}}$.

For example, the occurrence of $q$ in the source ${p\kon q}$ and
the occurrence of $q$ in the target ${q\kon p}$ of

\[
\c{\kon}_{p,q}\!\cirk(\Sd_{p,p}\kon \:\mj_q)\cirk(d_{p,\neg
p,p}\kon\:\mj_q)\cirk(\Dk_{p,p}\kon\:\mj_q)
\]

\noindent are directly linked in this headed factorized arrow term
of \PN, while the two occurrences of $p$ in its source and target
are not directly linked.

Take a headed factorized arrow term of \PN\ of the form ${g_2\cirk
f\cirk g_1}$ where $g_1$ is a $\Dk_{p,D}$-term and $g_2$ is
$\Sd_{p,C}$-term. Let ${\neg x_1\dis x_2}$ be the crown of the
head of $g_1$ (so $x_1$ and $x_2$ are both occurrences of $p$) and
let ${y_2\kon\neg y_1}$ be the crown of the head of $g_2$ (so
$y_1$ and $y_2$ are also occurrences of the same letter $p$). We
say that $g_1$ and $g_2$ are {\it confronted}\index{confronted
arrow terms} through $f$ when $x_i$ and $y_i$ are directly linked
for some ${i\in\{1,2\}}$ in the arrow term $f$.

Let a $\Dk_{p,A}$-term that is a factor of a factorized arrow term
$f$ be called a $\Dk$-{\it
factor}.\index{Deltaconfactor@$\Dk$-factor} We have an analogous
definition of $\Sd$-{\it
factor}\index{Sigmadis-factor@$\Sd$-factor} obtained by replacing
$\Dk$ by $\Sd$. We can then prove the following lemma.

\prop{Confrontation Lemma}{For\index{Confrontation Lemma} every
headed factorized arrow term ${g_2\cirk f\cirk g_1}$ of \PN\ such
that $g_1$ and $g_2$ are confronted through $f$ there is a headed
factorized arrow term $h$ of \PN\ with a subterm of the form
${g_2'\cirk f'\cirk g_1'}$ such that $g_1'$ is a $\Dk$-factor,
$g_2'$ is a $\Sd$-factor, $g_1'$ and $g_2'$ are confronted through
$f'$, and, moreover,}

\vspace{-1em}

\nav{(1)}{\it $f'$ is an arrow term of $\mds^{\neg p}\!$,}

\vspace{-1ex}

\nav{(2)}{\it ${g_2\cirk f\cirk g_1=h}$ in \PN,}

\vspace{-1ex}

\nav{(3)}{\it the number of $\Dk$-factors is equal in ${g_2\cirk
f\cirk g_1}$ and $h$, and the same for the number of
$\Sd$-factors.}

\vspace{2ex}

\dkz We proceed by induction on the number $n$ of factors of $f$
that are $\Dk$-factors or $\Sd$-factors. If $n=0$, then the arrow
term $f'$ coincides with the arrow term $f$.

If $n>0$, then let ${g_2\cirk f\cirk g_1}$ be of the form
${f_2\cirk g\cirk f_1}$ for $g$ a $\Dk_{q,E}$-term (we proceed
analogously when $g$ is a $\Sd_{q,E}$-term). According to the
stem-increasing equations of \S 6, we may assume that $g$
coincides with its head $\Dk_{q,E}$. Then by $\mbox{($\Dk$ {\it
nat})}$ we obtain in \PN\

\[
g_2\cirk f\cirk g_1=f_2\cirk(f_1\kon\mj_{\neg q\dis
q})\cirk\Dk_{q,E'}\!.
\]

\noindent After ${f_1\kon\mj_{\neg q\dis q}}$ in
${f_2\cirk(f_1\kon\mj_{\neg q\dis q})}$ is replaced by a headed
factorized arrow term ${g_2\cirk f''\cirk(g_1\kon\mj_{\neg q\dis
q})}$, we may apply the induction hypothesis to this arrow term,
because it can easily be seen that ${g_1\kon\mj_{\neg q\dis q}}$
and $g_2$ are confronted through $f''$, and $f''$ has one
$\Dk$-factor less than $f$. \qed

\vspace{2ex}

A headed factorized arrow term of \PN\ that has no subterm of the
form ${g_2\cirk f\cirk g_1}$ with $g_1$ and $g_2$ confronted
through $f$ is called {\it pure}.\index{pure arrow term} For a
pure arrow term $f$ there is a one-to-one correspondence, which we
call the \mbox{$\Dk$-{\it cap bijection}},\index{Deltaconcap
bijection@$\Dk$-cap bijection} between the $\Dk$-factors of $f$
and the caps of the partition part($Gf$). In this bijection, a cap
ties, in an obvious sense, the occurrences of $p$ in the crown
$\neg p\dis p$ of the head of the corresponding $\Dk$-factor.
There is an analogous one-to-one correspondence, which we call the
$\Sd$-{\it cup bijection},\index{Sigmadis-cup bijection@$\Sd$-cup
bijection} between the $\Sd$-factors of $f$ and the cups of
part($Gf$) (see \S 4 for the notions of cup and cap). Intuitively
speaking, this follows from the fact that in a sequence of cups
and caps tied to each other as in the following example:

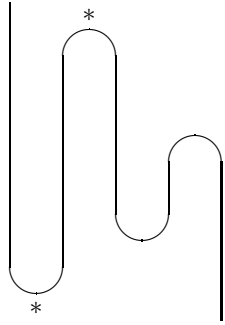
\begin{center}
\begin{picture}(80,120)

\put(0,20){\line(0,1){100}} \put(20,20){\line(0,1){80}}
\put(40,40){\line(0,1){60}} \put(60,40){\line(0,1){20}}
\put(80,0){\line(0,1){60}}

\put(10,20){\oval(20,20)[b]} \put(50,40){\oval(20,20)[b]}
\put(30,100){\oval(20,20)[t]} \put(70,60){\oval(20,20)[t]}

\put(10,7){\makebox(0,0)[t]{$\ast$}}
\put(30,113){\makebox(0,0)[b]{$\ast$}}

\end{picture}
\end{center}

\noindent cups and caps must alternate. For a pair made of a cap
and a cup that is its immediate neighbour, like those marked with
$\ast$ in the picture, we can find a subterm ${g_2\cirk f\cirk
g_1}$ such that $g_1$ and $g_2$ are confronted through $f$.

We can then prove the following.

\prop{Purification Lemma}{Every arrow term of \PN\ is equal in
\PN\ to a pure arrow term of \PN.}\index{Purification Lemma}

\dkz We apply first the Development Lemma of \S 3. If in the
resulting developed arrow term $h$ we have a subterm ${g_2\cirk
f\cirk g_1}$ with $g_1$ and $g_2$ confronted through $f$, then we
apply first the Confrontation Lemma to obtain a developed arrow
term $h'$ with a subterm of the form ${g_2'\cirk f'\cirk g_1'}$
where $g_1'$ and $g_2'$ are confronted through $f'$, and $f'$ is
an arrow term of $\mds^{\neg p}$.

Suppose that ${\neg x_2\dis x_3}$ is the crown of the head of
$g_1'$, and ${y_1\kon\neg y_2}$ is the crown of the head of
$g_2'$. Suppose $x_2$ is linked to $y_2$ in $f'$. Then, by Lemma~3
of \S 5, it is impossible that $x_3$ is linked to $y_1$, and so
there must be an occurrence $x_1$ of $p$ different from $x_3$ in
the source of $f'$ such that $x_1$ is linked to $y_1$ in $f'$, and
there must be an occurrence $y_3$ of $p$ different from $y_1$ in
the target of $f'$ such that $x_3$ is linked to $y_3$ in $f'$.
Next we apply the \emph{p-$\neg$p-p} Lemma of \S 6 to conclude
that ${g_2'\cirk f'\cirk g_1'}$ is equal to an arrow term $h''$ of
$\mds^{\neg p}$.

After replacing ${g_2'\cirk f'\cirk g_1'}$ in $h'$ by $h''$, we
obtain a headed factorized arrow term in which there is one
$\Dk$-factor and one $\Sd$-factor less than in $h'$, and hence
also than in $h$, by clause (3) of the Confrontation Lemma.

If $x_3$ is linked to $y_1$, then we reason analogously by
applying Lemma~3 of \S 5 and the \emph{$\neg$p-p-$\neg$p} Lemma of
\S 6.

We can iterate this procedure, which must terminate, because the
number of $\Dk$-factors and $\Sd$-factors in $h$ is finite. \qed

\vspace{2ex}

We can then prove the following.

\prop{\PN\ Coherence}{The functor $G$ from \PN\ to Br is
faithful.}\index{PN Coherence@\PN\ Coherence}

\dkz Suppose for $f$ and $g$ arrow terms of \PN\ of the same type
${A\vdash B}$ we have ${Gf=Gg}$. By the Purification Lemma, we can
assume that $f$ and $g$ are pure arrow terms. Since ${Gf=Gg}$, by
the $\Dk$-cap and $\Sd$-cup bijections we must have the same
number ${n\geq 0}$ of $\Dk$-factors in $f$ and $g$ and the same
number ${m\geq 0}$ of $\Sd$-factors in $f$ and $g$. We proceed by
induction on ${n\pl m}$.

If ${n\pl m=0}$, then we just apply \mds\ Coherence. Suppose now
${n>0}$. So there is a $\Dk$-factor in $f$ and a $\Dk$-factor in
$g$ that correspond by the $\Dk$-cap bijections to the same cap of
part($Gf$), which is equal to part($Gg$). By using the head
increasing equations of \S 6, together with $\mbox{($\Dk$ {\it
nat})}$, we obtain in \PN\

\[
f=f'\cirk\!\Dk_{p,A},\quad\quad\quad  g=g'\cirk\!\Dk_{p,A}
\]

\noindent for $f'$ and $g'$ pure arrow terms of the same type
${A\kon(\neg p\dis p)\vdash B}$, and such that the number of
$\Dk$-factors in $f'$ and $g'$ is ${n\mn 1}$ in each, and the
number of $\Sd$-factors in $f'$ and $g'$ is $m$ in each, the same
number we had for the $\Sd$-factors in $f$ and $g$. So we have

\[
G(f'\cirk\!\Dk_{p,A})=Gf=Gg=G(g'\cirk\!\Dk_{p,A}).
\]

We can show that ${Gf'=Gg'}$. We obtain $Gf'$ out of
${G(f'\cirk\!\Dk_{p,A})}$ in the following manner. We first remove
from the partition part(${G(f'\cirk\!\Dk_{p,A})}$) a cap
${\{k_t,l_t\}}$, where the \mbox{$k\pl 1$-th} occurrence of letter
in $B$ is an occurrence of $p$ in a subformula $\neg p$ of $B$,
and the \mbox{$l\pl 1$-th} occurrence of letter in $B$ is an
occurrence of $p$ that is not in a subformula $\neg p$ of $B$
(here we have either ${k<l}$ or ${l<k}$). After this removal, we
add two new transversals:

\[
\{GA_s,k_t\},\quad \{(GA\pl 1)_s,l_t\},
\]

\noindent and this yields part($Gf'$). Since $Gg'$ is obtained
from ${G(g'\cirk\!\Dk_{p,A})}$, which is equal to
${G(f'\cirk\!\Dk_{p,A})}$ in exactly the same manner, we obtain
that ${Gf'=Gg'}$.

Then, by the induction hypothesis, we have that ${f'=g'}$ in \PN,
which implies that ${f=g}$ in \PN. We proceed analogously in the
induction step when ${m>0}$, via $\Sd$-factors. \qed

\vspace{2ex}

From \PN\ Coherence and the equivalence between the categories
\PNN\ and \PN, proved in the preceding section, we may conclude in
the following manner that the functor $G$ from \PNN\ to \emph{Br}
is faithful.

\vspace{2ex}

\noindent{\sc Proof of \PNN\ Coherence.} Suppose that for $f$ and
$g$ arrows of \PNN\ of the same type we have ${Gf=Gg}$. Then, as
we noted at the beginning of this section, we have ${GFf=GFg}$,
and hence ${Ff=Fg}$ in \PN\ by \PN\ Coherence. It follows that
${f=g}$ in \PNN\ by the equivalence of the categories \PNN\ and
\PN. \qed

\vspace{2ex}

\noindent So we have proved \PNN\ Coherence, announced at the end
of \S 4.

\section{\large\bf The category \mmds}

In this and in the next section we add \emph{mix}\index{mix}
arrows of the type ${A\kon B\vdash A\dis B}$ to proof-net
categories, together with appropriate conditions that will enable
us to prove coherence with respect to \emph{Br} for the resulting
categories, which we call mix-proof-net categories. The mix
arrows, which underly the mix principle of linear logic, were
treated extensively in \cite{DP04} (Chapters 8, 10, 11, 13). The
proof of coherence for mix-proof-net categories is an adaptation
of the proof of coherence for proof-net categories given in the
preceding sections.

The category \mmds\ is defined as the category \mds\
in\index{MDS@\mmds\ category} \S 2 save that we have the
additional primitive arrow terms

\[
m_{A,B}\!:A\kon B\vdash A\dis B\index{m@$m$ arrow}
\]

\noindent for all objects, i.e.\ for all formulae, $A$ and $B$ of
$\eL_{\kon,\dis}$, and we assume the following additional
equations:

\begin{tabbing}
\quad\=($m$~{\it nat})\quad\=$m_{A\kon
B,C}\cirk\!\b{\kon}{\str}_{A,B,C}\;$\=$=m_{D,E}\cirk(f\kon
g)$\kill

\>($m$~{\it nat})\>$ (f\dis g)\cirk m_{A,B} = m_{D,E}\cirk(f\kon
g)$, \quad for $f\!:A\vdash D$ and $g\!:B\vdash
E$,\index{m1@($m$~{\it nat}) equation}
\\[1.5ex]
\>$(\b{\kon}{}\!m)$\>$m_{A\kon
B,C}\cirk\!\b{\kon}{\str}_{A,B,C}$\> $=d_{A,B,C}\cirk (\mj_A\kon
m_{B,C})$,\index{bconm@${(\b{\kon}{}m)}$ equation}
\\*[1ex]
\>$(\b{\dis}{}\!m)$ \> $\b{\dis}{\str}_{C,B,A}\!\cirk m_{C,B\dis
A}$\>$=(m_{C,B}\dis \mj_A)\cirk
d_{C,B,A}$,\index{bdism@${(\b{\dis}{}m)}$ equation}
\\*[1.5ex]
\>$(cm)$\>$m_{B,A}\cirk\!\c{\kon}_{A,B}\;=\;\c{\dis}_{B,A}\!\cirk
m_{A,B}$.\index{cm@${(cm)}$ equation}

\end{tabbing}

\noindent The proof-theoretical principle underlying $m_{A,B}$ is
called \emph{mix}\index{mix} (see \cite{DP04}, Section 8.1, and
references therein).

To obtain the functor $G$\index{G functor@$G$ functor} from \mmds\
to \emph{Br}, we extend the definition of the functor $G$ from
\mds\ to \emph{Br} (see \S 4) by adding the clause that says that
${Gm_{A,B}}$ is the identity arrow $\mj_{GA+GB}$ of \emph{Br}.
Then we have the following result of \cite{DP04} (Section 8.4).

\prop{\mmds\ Coherence}{The functor $G$ from $\mmds$ to Br is
faithful.}\index{MDS Coherence@\mmds\ Coherence}

In the remainder of this section we will prove some lemmata
concerning \mmds, which we will use for the proof of coherence in
the next section. For that we need some preliminaries.

For $x$ a particular proper subformula of a formula $A$ of
$\eL_{\kon,\dis}$, and $\!\ks\!\in\{\kon,\dis\}$, we define
$A^{-x}$ inductively as follows:

\begin{tabbing}

\quad\quad\quad\quad\quad\quad\quad\quad\quad\quad\=$(B\ks
x)^{-x}\;$\=$=(x\ks B)^{-x}=B$,
\\[2ex]
\quad for $x$ a proper subformula of $C$,
\\[1.5ex]
\>$(B\ks C)^{-x}$\>$=B\ks C^{-x}$,
\\[1ex]
\>$(C\ks B)^{-x}$\>$=C^{-x}\ks B$.

\end{tabbing}

For $i\in\{1,2\}$, let $A_i$ be a formula of $\eL_{\kon,\dis}$
with a proper subformula $x_i$, which is an occurrence of a letter
$q$, and let $x_i$ be the \mbox{$n_i$-th} occurrence of letter
counting from the left. We define the following functions
${\mu_i\!:\mbox{\boldmath $N$}-\{n_i\mn 1\}\str \mbox{\boldmath
$N$}}$:

\[
\mu_i(n)=_{\df}\left\{\begin{array}{ll} n & \mbox{\rm{if }} n<
n_i\mn 1
\\ n\mn 1 & \mbox{\rm{if }} n> n_i\mn 1.
\end{array} \right.
\]

The definition of \emph{linked} occurrence of a letter in an arrow
of \mmds\ is analogous to what we had in \S 5. Then we can prove
the following.

\prop{Lemma~1}{For every arrow term ${f\!:A_1\vdash A_2}$ of
\mmds\ such that $x_1$ and $x_2$ are linked in the arrow $f$,
there is an arrow term ${f^{-q}\!:A_1^{-x_1}\vdash A_2^{-x_2}}$ of
\mmds\ such that the members of \emph{part($Gf^{-q}$)} are
${\{s(\mu_1(m_1)),t(\mu_2(m_2))\}}$ for each ${\{s(m_1),t(m_2)\}}$
in \emph{part($Gf$)}, provided ${m_i\neq n_i\mn 1}$.}

\dkz We proceed by induction on the complexity of the arrow term
$f$. If $f$ is a primitive arrow term $\alpha_{B_1,\dots,B_m}$,
then for some ${j\in\{1,\ldots,m\}}$ we have that $x_i$ occurs in
a subformula $B_j$ of $A_i$. If $x_i$ is a proper subformula of
this subformula $B_j$, then $B^{-x_{i}}_j$ is defined, and
$f^{-q}$ is

\[
\alpha_{B_1,\dots,B_{j-1},B^{-x_{i}}_j,B_{j + 1},\dots,B_m}
\]

\noindent (note that $B^{-x_{1}}_j$ and $B^{-x_{2}}_j$ are the
same formula). If $x_i$ is not a proper subformula of the
subformula $B_j$, then  $f^{-q}$ is $\mj_{A^{-x_i}_i}$.

If $f$ is $g\cirk h$, then $f^{-q}$ is $g^{-q}\cirk h^{-q}$, and
if $f$ is $g\ks h$ for $\!\ks\!\in\{\kon,\dis\}$, then $f^{-q}$ is
either $g^{-q}\ks h$, or $g\ks h^{-q}$, or $g$ when $h$ is
$\mj_{x_1}$, or $h$ when $g$ is $\mj_{x_1}$. \qed

\vspace{2ex}

Here is an example of the application of Lemma~1. If
${f\!:A_1\vdash A_2}$ is

\begin{tabbing}
$((m_{q,p\kon
q}\cirk(\mj_q\:\kon\c{\kon}_{q,p})\cirk\!\c{\kon}_{q\kon
p,q})\dis\mj_p)\cirk d_{q\kon
p,q,p}\cirk\!\b{\kon}{\str}_{q,p,q\dis p}\cirk\!
\c{\kon}_{p\kon(q\dis p),q}:$
\\*[.5ex]
\`$(p\kon(q\dis p))\kon q\vdash (q\dis(p\kon q))\dis p$,
\end{tabbing}

\noindent where $x_1$ is the second (rightmost) occurrence of $q$
in $(p\kon(q\dis p))\kon q$, while $x_2$ is the second occurrence
of $q$ in $(q\dis(p\kon q))\dis p$, then ${f^{-q}\!:A_1^{-q}\vdash
A_2^{-q}}$ is

\begin{tabbing}
$((m_{q,p}\cirk(\mj_q\kon\mj_p)\cirk\!\c{\kon}_{p,q})\dis\mj_p)\cirk
d_{p,q,p}\cirk\mj_{p\kon(q\dis p)}\cirk\mj_{p\kon(q\dis p)}:$
\\*[.5ex]
\`$p\kon(q\dis p)\vdash (q\dis p)\dis p$,
\end{tabbing}

\noindent which is equal to
${((m_{q,p}\cirk\!\c{\kon}_{p,q})\dis\mj_p)\cirk d_{p,q,p}}$.

We define inductively a notion we call a \emph{context}:

\begin{itemize}
\item[] $\koc$ is a context; \item[] if $Z$ is a context and $A$ a
formula of $\eL_{\kon,\dis}$, then $Z\ks A$ and $A\ks Z$ are
contexts for ${\!\ks\!\in\{\kon,\dis\}}$.
\end{itemize}

Next we define inductively what it means for a context $Z$ to be
applied to an object $B$ of \mmds, which we write $Z(B)$, or to an
arrow term $f$ of \mmds, which we write $Z(f)$:

\begin{tabbing}
\mbox{\hspace{2em}}\= $(Z\ks A)(B)$ \= $=$ \=
\mbox{\hspace{10em}}\= $(Z\ks A)(f)$ \= $=$ \=\kill \> \>
$\koc(B)$\' $=$\> $B$,\> \> $\koc(f)$\' $=$\> $f$,
\\[1ex]
\> $(Z\ks A)(B)$\> $=$\> $Z(B)\ks A$,\> $(Z\ks A)(f)$\> $=$\>
$Z(f)\ks \mj_A$,
\\[1ex]
\> $(A\ks Z)(B)$\> $=$\> $A\ks Z(B)$;\> $(A\ks Z)(f)$\> $=$\>
$\mj_A\ks Z(f)$.
\end{tabbing}

\noindent We use $X$, $Y$, $Z,\dots$ for contexts.\index{X
context@$X$ context}\index{Y context@$Y$ context}\index{Z
context@$Z$ context}

For ${f\!:A\vdash C}$ an arrow of \mmds, we say that an occurrence
$x$ of a formula $B$ as a subformula of $A$ and an occurrence $y$
of the same formula $B$ as a subformula of $C$ are \emph{linked}
in $f$ when the \mbox{$n$-th} letter in $x$ is linked in $f$ to
the \mbox{$n$-th} letter in $y$.

Let ${f\!:X(p)\kon B\vdash Y(p\kon B)}$ be an arrow term of \mmds\
such that the displayed occurrences of $p$ in the source and
target, and also the displayed occurrences of $B$, are linked in
the arrow $f$. Then, by successive applications of Lemma~1, for
each occurrence of a letter in $B$, we obtain the arrow term
${f^{-B}\!:X(p)\vdash Y(p)}$ of \mmds, and the displayed
occurrences of $p$ in $X(p)$ and $Y(p)$ are linked in the arrow
${f^{-B}}$.

Let ${f^{\dag}\!:X(p\kon B)\vdash Y(p\kon B)}$ be the arrow term
of \mmds\ obtained from ${f^{-B}}$ by replacing the occurrences of
$p$ that correspond to those displayed in ${X(p)}$ and ${Y(p)}$ by
occurrences of ${p\kon B}$. This replacement is made in the
indices of primitive arrow terms that occur in ${f^{-B}}$, and it
need not involve all the occurrences of $p$ in these indices. For
example, if $X$ is ${\koc\kon(q\dis p)}$ and $Y$ is
${(q\dis\koc)\dis p}$, while ${f^{-B}}$ is

\[
((m_{q,p}\cirk\!\c{\kon}_{p,q})\dis\mj_p)\cirk
d_{p,q,p}\!:p\kon(q\dis p)\vdash (q\dis p)\dis p,
\]

\noindent then ${f^{\dag}}$ is

\[
((m_{q,p\kon B}\cirk\!\c{\kon}_{p\kon B,q})\dis\mj_p)\cirk
d_{p\kon B,q,p}\!:(p\kon B)\kon(q\dis p)\vdash (q\dis (p\kon
B))\dis p.
\]

Then we can prove the following.

\prop{Lemma~2$\kon$}{Let ${f\!:X(p)\kon B\vdash Y(p\kon B)}$ and
${f^{\dag}\!:X(p\kon B)\vdash Y(p\kon B)}$ be as above. Then there
is an arrow term ${h_X\!:X(p)\kon B\vdash X(p\kon B)}$ of \mds\
such that ${f=f^{\dag}\cirk h_X}$ in \mmds.}

\dkz We construct the arrow term $h_X$\index{hax@$h_X$ arrow} of
\mds\ by induction on the complexity of the context $X$. For the
basis we have that $h_{\koci}$ is $\mj_{p\kon B}$. In the
induction step we have

\begin{tabbing}
\quad\quad\quad\quad\=$h_{Z\kon
A}\;$\=$=(h_Z\kon\mj_A)\cirk\!\c{\kon}_{A,Z(p)\kon
B}\!\cirk\!\b{\kon}{\rts}_{A,Z(p),B}\!\cirk(\c{\kon}_{Z(p),A}\kon\:\mj_B)$,
\\*[.5ex]
\>$h_{Z\dis A}$\>$=(h_Z\dis\mj_A)\cirk\!\c{\dis}_{Z(p)\kon
B,A}\!\cirk \,d^R_{A,Z(p),B}\cirk(\c{\dis}_{A,Z(p)}\kon\:\mj_B)$,
\\[.5ex]
\>$h_{A\kon Z}$\>$=(\mj_A\kon
h_Z)\cirk\!\b{\kon}{\rts}_{A,Z(p),B}$,
\\*[.5ex]
\>$h_{A\dis Z}$\>$=(\mj_A\dis h_Z)\cirk \,d^R_{A,Z(p),B}$.
\end{tabbing}

It is easy to see that ${Gf=G(f^{\dag}\cirk h_X)}$, and then the
lemma follows by applying \mmds\ Coherence. \qed

\vspace{2ex}

Let ${f\!:Y(B\dis p)\vdash B\dis X(p)}$ be an arrow term of \mmds\
such that the displayed occurrences of $p$ in the source and
target, and also the displayed occurrences of $B$, are linked in
the arrow $f$. Then, as above by Lemma~1, we obtain the arrow term
${f^{-B}\!:Y(p)\vdash X(p)}$ of \mmds, and the displayed
occurrences of $p$ in $Y(p)$ and $X(p)$ are linked in the arrow
${f^{-B}}$.

Let ${f^{\dag}\!:Y(B\dis p)\vdash X(B\dis p)}$ be the arrow term
of \mmds\ obtained from ${f^{-B}}$ by replacing the occurrences of
$p$ that correspond to those displayed in ${Y(p)}$ and ${X(p)}$ by
occurrences of ${B\dis p}$ (cf.\ the example above). Then we can
prove the following, analogously to Lemma~2$\kon$.

\prop{Lemma~2$\dis$}{Let ${f\!:Y(B\dis p)\vdash B\dis X(p)}$ and
${f^{\dag}\!:Y(B\dis p)\vdash X(B\dis p)}$ be as above. Then there
is an arrow term ${h_X\!:X(B\dis p)\vdash B\dis X(p)}$ of \mds\
such that ${f=h_X\cirk f^{\dag}}$ in \mmds.}

\section{\large\bf \MPNN\ Coherence}

The category \MPNN\ is\index{MPNN@\MPNN\ category} defined as the
category \PNN\ in \S 3 save that we have the additional primitive
arrow terms ${m_{A,B}\!:A\kon B\vdash A\dis B}$ for all objects
$A$ and $B$ of \PNN, and we assume as additional equations
($m$~{\it nat}), (${\b{\kon}{}\!m}$), (${\b{\dis}{}\!m}$) and
(${cm}$) of the preceding section. To obtain the functor
$G$\index{G functor@$G$ functor} from \MPNN\ to \emph{Br}, we
extend the definition of the functor $G$ from \PNN\ to \emph{Br}
by adding the clause that says that ${Gm_{A,B}}$ is the identity
arrow $\mj_{GA+GB}$ of \emph{Br}.

A \emph{mix-proof-net}\index{mix-proof-net category} category is
defined as a proof-net category (see \S 3) that has in addition a
natural transformation $m$ satisfying the equations
(${\b{\kon}{}\!m}$), (${\b{\dis}{}\!m}$) and (${cm}$). The
category \MPNN\ is up to isomorphism the free mix-proof-net
category generated by \Pe.

The category \MPN\ is\index{MPN@\MPN\ category} defined as the
category \PN\ in \S 6 save that we have the additional primitive
arrow terms ${m_{A,B}}$ for all objects of \PN, and we assume as
additional equations ($m$~{\it nat}), (${\b{\kon}{}\!m}$),
(${\b{\dis}{}\!m}$) and (${cm}$). We can prove that \MPNN\ and
\MPN\ are equivalent categories as in \S 7. (We have an additional
case involving ${m_{A,B}}$ in the proof of the analogue of the
Auxiliary Lemma of \S 7, and similar trivial additions elsewhere;
otherwise the proof is quite analogous.)

We have a functor $G$\index{G functor@$G$ functor} from \MPN\ to
\emph{Br} defined by restricting the definition of the functor $G$
from \MPNN\ to \emph{Br} (cf.\ the beginning of \S 8), and we will
prove the following.

\prop{\MPN\ Coherence}{The functor $G$ from \MPN\ to Br is
faithful.}\index{MPN Coherence@\MPN\ Coherence}

The proof of this coherence proceeds as the proof of \PN\
Coherence in \S 8. The only difference is in the $\Xk$-Permutation
and $\Xd$-Permutation Lemmata of \S 6.

The formulation of the $\Xk$-Permutation
Lemma\index{Xicon-Permutation Lemma@$\Xk$-Permutation Lemma} is
modified by replacing \PN\ and $\mds^{\neg p}$ by respectively
\MPN\ and $\mmds^{\neg p}$, where the category $\mmds^{\neg p}$ is
defined as \mmds\ save that it is generated not by \Pe, but by
${\Pe\cup\Pe^{\neg}}$ (cf.\ \S 6); moreover, we assume that $y_1$
and $\neg y_2$ occur in $E$ within a subformula of the form
${p\kon(\neg y_2\dis y_1)}$ or ${\neg p\kon(y_1\dis\neg y_2)}$. We
modify the proof of this lemma as follows.

If in $E$ we have ${p\kon(\neg y_2\dis y_1)}$, then by the
stem-increasing equations of \S 6 we have that the
$\Xk_{p,B}$-term ${g\!:C\vdash D}$ is equal to
${f''\cirk\!\Dk_{p,C}}$ for ${f''\!:C\kon(\neg p\dis p)\vdash D}$
an arrow term of $\mds^{\neg p}$, and so for ${f\!:D\vdash E}$ an
arrow term of $\mmds^{\neg p}$ satisfying the conditions of the
lemma we have in $\MPN$

\[
f\cirk g=f\cirk f''\cirk\!\Dk_{p,C}\!.
\]

Then we apply Lemma~2$\kon$ of the preceding section to

\[
f\cirk f''\!:C\kon(\neg p\dis p)\vdash E,
\]

\noindent where $C$ is $X(p)$, ${\neg p\dis p}$ is $B$ and $E$ is
${Y(p\kon(\neg p\dis p))}$. So for

\[
h_X\!:X(p)\kon(\neg p\dis p)\vdash X(p\kon(\neg p\dis p))
\]

\noindent an arrow term of $\mds^{\neg p}$, and

\[
(f\cirk f'')^{\dag}\!:X(p\kon(\neg p\dis p))\vdash Y(p\kon(\neg
p\dis p))
\]

\noindent we have

\[
f\cirk f''=(f\cirk f'')^{\dag}\cirk h_X.
\]

By the $\Xk$-Permutation Lemma of \S 6 we have

\[
h_X\cirk\!\Dk_{p,C}\;=g'\cirk f'
\]

\noindent where $g'$ is the $\Dk_{p,p}$-term $X(\Dk_{p,p})$, and
by bifunctorial and naturality equations we have

\[
(f\cirk f'')^{\dag}\cirk X(\Dk_{p,p})=Y(\Dk_{p,p})\cirk (f\cirk
f'')^{-(\neg p\dis p)}.
\]

\noindent Note that ${(f\cirk f'')^{\dag}}$ is obtained from
${(f\cirk f'')^{-(\neg p\dis p)}\!:X(p)\vdash Y(p)}$ by
replacement of $p$.

So we have in $\MPN$

\[
\begin{array}{ll}
 f\cirk g \!\!\!& = f\cirk f''\cirk\!\Dk_{p,C}\\[.5ex]
 &=(f\cirk f'')^{\dag}\cirk h_X\cirk\!\Dk_{p,C}\\[.5ex]
 &=(f\cirk f'')^{\dag}\cirk X(\Dk_{p,p})\cirk f'\\[.5ex]
 &=Y(\Dk_{p,p})\cirk f'''
\end{array}
\]

\noindent for $f'''$, which is ${(f\cirk f'')^{-(\neg p\dis
p)}\cirk f'}$, an arrow term of $\mmds^{\neg p}$.

We proceed analogously if in $E$ we have ${\neg p\kon(y_1\dis\neg
y_2)}$; instead of $\Dk_{p,p}$ we then have $\Dkp_{p,p}$. We have
an analogous reformulation of the $\Xd$-Permutation
Lemma\index{Xidis-Permutation Lemma@$\Xd$-Permutation Lemma} of \S
6, with a proof based on Lemma~2$\dis$ of the preceding section.

Instead of Lemma~2$\kon$ of the preceding section, we could have
proved, with more difficulty, an analogous lemma where $f$ is of
type

\[
Z(X_1(p)\kon X_2(B))\vdash Y(p\kon B),
\]

\noindent and $f^{\dag}$ is of one of the following types:

\[
\begin{array}{cc} Z(X_1(p\kon B)\kon (X_2(B))^{-B})\vdash Y(p\kon
B),
\\*[1.5ex]
Z(X_1(p\kon B))\vdash Y(p\kon B).
\end{array}
\]

\noindent Then in the proof of the $\Xk$-Permutation Lemma
modified for \MPN\ we would not need to pass from $g$ to
${f''\cirk\!\Dk_{p,C}}$ via stem-increasing equations, but this
alternative approach is altogether less clear.

Note that we have no analogue of Lemma~2 of \S 5 for \mmds. The
lack of this lemma, on which we relied in \S 6 for the proof of
the $\Xk$-Permutation and $\Xd$-Permutation Lemmata, is tied to
the modifications we made for these lemmata with \MPN. We have
also no analogue of Lemma~4 of \S 5, but the analogue of Lemma~3
of \S 5 does hold.

From \MPN\ Coherence and the equivalence of the categories \MPNN\
and \MPN\ we can then infer the following.

\prop{\MPNN\ Coherence}{The functor $G$ from \MPNN\ to Br is
faithful.}\index{MPNN Coherence@\MPNN\ Coherence}

\end{document}